# Doubles of groups and hyperbolic LERF 3-manifolds

By Rita Gitik*

### Abstract

We show that the quasiconvex subgroups in doubles of certain negatively curved groups are closed in the profinite topology. This allows us to construct the first known large family of hyperbolic 3-manifolds such that any finitely generated subgroup of the fundamental group of any member of the family is closed in the profinite topology.

### Introduction

The profinite topology on a group $G$ is defined by proclaiming all finite index subgroups of $G$ to be the base open neighborhoods of the identity in $G$. We denote it by $PT(G)$. A group $G$ is RF (residually finite) if the trivial subgroup is closed in $PT(G)$, which happens if and only if $PT(G)$ is Hausdorff. A group $G$ is LERF (locally extended residually finite) if any finitely generated subgroup of $G$ is closed in $PT(G)$. RF and LERF groups have been studied for a long time, and they have various important properties. For example, finitely generated RF groups have solvable word problem and finitely generated LERF groups have solvable generalized word problem; see [A-G], [B-B-S], [Gi 2] and [We] for various results and additional references. The class of RF groups is very rich. It contains all finitely generated linear groups and all fundamental groups of geometric 3-manifolds. However, few examples of LERF groups were known.

We say that a 3-manifold is LERF if its fundamental group is LERF, and we say that a 3-manifold with boundary is hyperbolic if its interior has a complete hyperbolic structure. In this paper we construct the first known large nontrivial class of hyperbolic LERF 3-manifolds with boundary, and a new large class of closed hyperbolic 3-manifolds, which have all their surface

*Research partially supported by NSF grant DMS 9022140 at MSRI.

subgroups and all their geometrically finite subgroups closed in the profinite topology.

If the fundamental group of a compact orientable irreducible 3-manifold $M$ has a surface subgroup $S$ which is closed in the profinite topology of $\pi_1(M)$, then $M$ is virtually Haken. Specifically, there exists a finite cover $N$ of $M$ such that $S$ is contained in $\pi_1(N)$ and is carried by a surface embedded in $N$. A conjecture of Waldhausen asserts that any such closed 3-manifold $M$ whose fundamental group contains a surface subgroup is virtually Haken, hence the importance of the LERF property in 3-manifolds.

It was conjectured that all finitely generated 3-manifold groups are LERF, and P. Scott proved in [Sco 1, 2] that compact Seifert fibered spaces are LERF. However, a non-LERF compact graph manifold was described in [B-K-S], and it appears that most graph manifolds are not LERF, ([L-N], [R-W]). Still, little was known about hyperbolic LERF 3-manifolds. M. Hall proved in [Hall] that free groups are LERF, so that handlebodies are LERF. P. Scott proved in [Sco 1] that surface groups are LERF, so that $I$-bundles over surfaces are LERF. He also showed that all geometrically finite subgroups of certain closed hyperbolic 3-manifolds are closed in the profinite topology. This limited information about the profinite topology on the fundamental groups of hyperbolic 3-manifolds prompted W. Thurston to ask in [Thu] whether finitely generated Kleinian groups are LERF or whether they have special subgroups closed in the profinite topology.

Since then it was shown in [B-B-S] that a free product of two free groups with cyclic amalgamation is LERF, so an annulus sum of two handlebodies is LERF. Later the author showed in [Gi 2] that the free product of a LERF group and a free group amalgamated over a cyclic group maximal in the free factor is LERF; hence the sum of any LERF hyperbolic 3-manifold and a handlebody along an annulus maximal in the handlebody is LERF.

The following theorem is the main topological result of this paper.

THEOREM 1. *Let $M$ be a compact hyperbolic LERF 3-manifold with boundary, which does not have boundary tori, let $B$ be a connected submanifold of the boundary of $M$, such that $B$ is incompressible in $M$, and let $D(M)$ be the double of $M$ along $B$. If $D(M)$ is hyperbolic, has nonempty boundary, and has no boundary tori, then $D(M)$ is LERF. If the boundary of $D(M)$ is empty, then any geometrically finite subgroup and any freely indecomposable geometrically infinite subgroup (hence any closed surface subgroup) of the fundamental group of $D(M)$ is closed in the profinite topology.*

Theorem 1 is a corollary of Theorem 2, and its proof is given at the end of Section 1. The "no boundary tori" condition seems not to be essential, and the author plans to remove it, at least in some cases, in a subsequent paper.

Theorem 1 enables us to construct hyperbolic 3-manifolds with LERF fundamental group as follows. Let $M$ be as in Theorem 1. Initial examples are handlebodies or $I$-bundles over closed surfaces of negative Euler characteristic, or annulus sums of several handlebodies with such an $I$-bundle. In general, the boundary of $M$ might be compressible (for example, if $M$ is a handlebody) or $M$ might be not acylindrical (for example, if $M$ is an $I$-bundle over a closed surface). If $M$ has incompressible boundary and is not acylindrical, we can use the characteristic submanifold theorem of Jaco-Shalen and Johannson to show that any boundary component of $M$ carries many essential simple closed curves $C$ which separate this boundary component in two parts $A$ and $B$, each incompressible in $M$, such that $\pi_1(A)$ and $\pi_1(B)$ are malnormal subgroups of $\pi_1(M)$. Then Theorem 1 implies that the double of $M$ along either $A$ or $B$ is LERF. As $D(M)$ has nonempty boundary, we can apply the characteristic submanifold theorem to a boundary component of $D(M)$, and double $D(M)$ along a part of its boundary, creating a hyperbolic LERF manifold $D(D(M))$. Iteration of this process produces a large family of hyperbolic LERF 3-manifolds with boundary.

In order to construct a closed hyperbolic manifold $N$ such that any geometrically finite subgroup of $\pi_1(N)$ is closed in $\operatorname{PT}(\pi_1(N))$, we need to start with $M$, as in Theorem 1, such that its boundary is connected and incompressible. If the boundary of $M$ is acylindrical (for example, totally geodesic), then the double of $M$ along the whole boundary will be hyperbolic and closed, hence it will have the required properties.

If the boundary of $M$ is not acylindrical, we still can carry the construction, but in two steps. We need to find a simple closed essential curve $C$ separating the boundary of $M$ in two parts $A$ and $B$ satisfying much stricter conditions, namely:

1) $\pi_1(A)$ is a malnormal subgroup of $\pi_1(M)$.

2) $\pi_1(D(B))$ is a malnormal subgroup of $\pi_1(D(M))$, where $D(M)$ is the double of $M$ over $A$, and $D(B)$ is the double of $B$ over $C$.

Then the double $N$ of $D(M)$ over $D(B)$ is a closed hyperbolic 3-manifold with the required properties.

We can take $M$ to be a twisted $I$-bundle over a nonorientable surface of genus 2, because there exist separating curves $C$ in its boundary such that the groups $\pi_1(A)$ and $\pi_1(B)$ inject in $\pi_1(M)$, $M$ has no essential cylinders with both ends in $A$, $M$ has no essential cylinders with both ends in $B$ and $M$ has no cylinders connecting $A$ and $B$, hence $A$ and $B$ have properties 1) and 2) mentioned above.

## 1. The profinite topology on doubles of groups

The main group-theoretical result of this paper is a combination theorem for the profinite topology on a special class of groups. It is well-known that free products preserve RF and LERF groups, but free products with amalgamation usually do not (cf. [A-G], [L-N]). It is shown in [G-R 1] that adjunction of roots need not preserve the property LERF, so one should not expect the profinite topology on groups to behave reasonably even under free products with cyclic amalgamation. In this paper we study the profinite topology on a special class of amalgamated free products, called doubles. The graph-theoretical techniques developed in this paper and in [Gi 2] allow the author to prove new combination theorems (not only about doubles) for profinite topology on groups. As these results are not connected with the main subject of this paper, they will be described somewhere else.

*Definition* 1.1. Let $G_0$ be a subgroup of a group $G$, let $H$ be an isomorphic copy of $G$ with a fixed isomorphism $\alpha : G \to H$ and let $H_0 = \alpha(G_0)$. The double of $G$ along $G_0$ is the amalgamated free product $D = G *_{G_0 = H_0} H$. We call $G$ and $H$ "the factors of $D$". When $X$ is a generating set of $G$, then $Y = \alpha(X)$ is a generating set of $H$.

The following example shows that a subgroup of $G$ which is closed in $\mathrm{PT}(G)$ does not have to be closed in $\mathrm{PT}(D)$.

*Example* 1.2. A double of an RF group need not be RF. Let $G = \langle a, c | a^{-1}cac^{-2} \rangle$ and let $G_0 = \langle c \rangle$. The group $G$ is RF, but it is shown in [Hi] that the double $D$ of $G$ along $G_0$ is not. Hence the trivial subgroup is closed in $\mathrm{PT}(G)$, but it is not closed in $\mathrm{PT}(D)$. Note that $G_0$ is not closed in $\mathrm{PT}(G)$, because the element $aca^{-1}$ belongs to the closure of $G_0$ in $\mathrm{PT}(G)$.

This example is generic, as D. Long and G. Niblo proved in [L-N] that the double of an RF group $G$ along $G_0$ is RF if and only if $G_0$ is closed in $\mathrm{PT}(G)$. The following more general statement is proved in [Gi 5].

LEMMA 1.3. *Let $D$ be the double of $G$ along $G_0$. If $G_0$ is closed in $\mathrm{PT}(G)$, then any subgroup of $G$ which is closed in $\mathrm{PT}(G)$ is closed in $\mathrm{PT}(D)$. If $G_0$ is not closed in $\mathrm{PT}(G)$, then no subgroup of $G$ is closed in $\mathrm{PT}(D)$.*

An obvious necessary condition for a subgroup $S$ of $D$ to be closed in $\mathrm{PT}(D)$ is that the intersection of $S$ with any conjugate of a factor of $D$ must be closed in the profinite topology of the conjugate. If $G$ is LERF, this condition holds if the intersection of $S$ with any conjugate of a factor of $D$ is finitely generated or, equivalently, the intersection of $S$ with any conjugate of $G_0$ is finitely generated. Of course, there exist infinitely generated subgroups which are closed in the profinite topology; however, detecting such subgroups seems to be a very difficult problem.

*Example* 1.4. A double of a LERF group need not be LERF. Let $F_n$ denote the free group of rank $n$. Let $G = F_1 \times F_2 = \langle u \rangle \times \langle x, y \rangle$, and let $G_0 = F_2 = \langle x, y \rangle$. It is shown in [A-G] that $G$ is LERF, but the double of $G$ along $G_0$, which is isomorphic to $F_2 \times F_2$, is not LERF, although it is RF.

Recall that a group $D$ has fgip (finitely generated intersection property) if the intersection of any pair of its finitely generated subgroups is finitely generated, and a subgroup $G_0$ of $D$ has fgip in $D$ if the intersection of $G_0$ with any finitely generated subgroup of $D$ is finitely generated. It is easy to exhibit a finitely generated subgroup of $F_2 \times F_2$ in Example 1.4 such that its intersection with the amalgamating subgroup $G_0$ is infinitely generated; hence the failure of $F_2 \times F_2$ to be LERF can be attributed to the failure of the amalgamating subgroup $G_0$ to have fgip in $F_2 \times F_2$. However, the situation is much more complicated, because there exists a double $D$ of $F_2$ along a finite index subgroup of $F_2$ such that $D$ has a subgroup isomorphic to $F_2 \times F_2$. Such a $D$ cannot be LERF (cf. [Ge], [Rips]). As a finite index subgroup has fgip in any finitely generated group, the problem can be caused only by the way the amalgamating subgroup $G_0$ is embedded in $G$.

In this paper we give a condition on $G_0$ which forces $D$ to be LERF. The main technical group-theoretical results of this paper are the following theorems.

THEOREM 4.4. *Let $S$ be a finitely generated subgroup of the double $D$ of a LERF group $G$ along a finitely generated subgroup $G_0$, such that the intersection of $S$ with any conjugate of $G_0$ is finitely generated. If $G_0$ is strongly separable (see Definition 4.2) in $G$, then $S$ is closed in $\mathrm{PT}(D)$. Hence if $G_0$ is strongly separable in $G$ and has fgip in $D$, then $D$ is LERF.*

THEOREM 5.4. *A finitely generated malnormal subgroup of a locally quasiconvex LERF negatively curved group is strongly separable.*

Recall that a group is locally quasiconvex if all its finitely generated subgroups are quasiconvex, and a subgroup $H$ is malnormal in $G$ if for any $g \notin H$ the intersection of $H$ and $gHg^{-1}$ is trivial.

Theorem 4.4 and Theorem 5.4 imply our main group-theoretical result.

THEOREM 2. *Let $G$ be a finitely generated locally quasiconvex negatively curved LERF group, and let $D$ be the double of $G$ along a finitely generated subgroup $G_0$. If $G_0$ is malnormal in $G$, then any quasiconvex subgroup of $D$ is closed in $\mathrm{PT}(D)$. Hence if $D$ is locally quasiconvex, then $D$ is LERF.*

*Proof.* Let $G$ be a finitely generated locally quasiconvex negatively curved group, and let $G_0$ be a finitely generated malnormal subgroup of $G$. Theorem 5.4 implies that $G_0$ is strongly separable in $G$. As $G_0$ is finitely generated, it is quasiconvex in $G$; hence $D$ is negatively curved ([B-F], [Gi 6]). Then it is shown in [Gi 3] that all conjugates of $G_0$ in $D$ are quasiconvex in $D$.

Let $S$ be a quasiconvex subgroup of $D$. As quasiconvex subgroups of finitely generated groups are finitely generated, and as the intersection of two quasiconvex subgroups is a quasiconvex subgroup ([Gre], [Gi 3]), it follows that the intersection of $S$ with any conjugate of $G_0$ is finitely generated. Therefore Theorem 4.4 implies that $S$ is closed in $\mathrm{PT}(D)$.

Theorem 2 easily implies Theorem 1, as follows.

*Proof of Theorem* 1. As $M$ is compact, its fundamental group is finitely generated. As $M$ is hyperbolic and has no boundary tori, its fundamental group is negatively curved. If $D(M)$ is hyperbolic, then $M$ does not contain essential cylinders with both ends in $B$, so $\pi_1(B)$ is a malnormal subgroup of $\pi_1(M)$. A theorem of W. Thurston states that if a hyperbolic 3-manifold with finitely generated fundamental group has at least one boundary component which is not a torus, then its fundamental group is locally quasiconvex. As $M$ has nonempty boundary and no boundary tori, $\pi_1(M)$ is locally quasiconvex.

If $D(M)$ has nonempty boundary and no boundary tori, then $\pi_1(D(M))$ is also locally quasiconvex and negatively curved; hence Theorem 2 implies that $D(M)$ is LERF.

If the boundary of $D(M)$ is empty, then Theorem 2 implies that any quasiconvex subgroup of $D(M)$ is closed in the profinite topology. A theorem of F. Bonahon ([Bo]) implies that any nonquasiconvex freely indecomposable subgroup of $\pi_1(D(M))$ is closed in the profinite topology. Hence any subgroup of $\pi_1(D(M))$ which is isomorphic to the fundamental group of a closed surface is closed in $\mathrm{PT}(\pi_1(D(M)))$.

It is shown in [Gi 3] that a double of a locally quasiconvex negatively curved group along a malnormal cyclic subgroup is locally quasiconvex. As fundamental groups of closed surfaces of genus greater than 1 are locally quasiconvex, negatively curved and LERF, the following statement is a special case of Theorem 2.

Corollary 1.5. *A double of a locally quasiconvex negatively curved* LERF *group (for example, a double of a fundamental group of a closed surface of genus greater than* 1) *along a malnormal cyclic subgroup is* LERF.

Note that there exist examples of non-LERF groups which are doubles of LERF groups over cyclic subgroups ([L-N], [G-R 1], [A-D]). As a cyclic subgroup has fgip in any group, this phenomenon is caused by the way the amalgamating subgroup $G_0$ is embedded in $G$. Niblo in [Ni] proved that if $D$ is a LERF group which is a double of a LERF group $G$ along $G_0$, then for any finitely generated subgroup $S$ of $G$ the set $G_0 S$ is closed in $\mathrm{PT}(G)$. He also showed that this condition on $G_0$ is not sufficient for $D$ to be LERF, even when $G_0$ is cyclic.

## 2. Preliminaries

This section contains a summary of graph-theoretical methods developed by the author in [Gi 1] and in [Gi 2]. The detailed proofs of the quoted results appeared in [Gi 2].

*Definition* 2.1.   Let $X$ be a set, let $X^* = \{x, x^{-1} | x \in X\}$, and for $x \in X$ define $(x^{-1})^{-1} = x$. Consider a group $G$ generated by the set $X$. Let $G_0$ be a subgroup of $G$, and let $\{G_0 g\}$ denote the set of right cosets of $G_0$ in $G$. The relative Cayley graph of $G$ with respect to $G_0$ (or the coset graph) is an oriented graph whose vertices are the right cosets $\{G_0 g\}$ and the set of edges is $\{G_0 g\} \times X^*$, such that an edge $(G_0 g, x)$ begins at the vertex $G_0 g$ and ends at the vertex $G_0 gx$. We denote it $\mathrm{Cayley}(G, G_0)$. Note that $G_0$ acts on the Cayley graph of $G$ by left multiplication, and $\mathrm{Cayley}(G, G_0)$ can be defined as the quotient of the Cayley graph of $G$ by this action.

Let $K$ be the standard 2-complex representing the group $G = \langle X | R \rangle$, i.e. $K$ has one vertex, $|X|$ oriented edges and $|R|$ 2-cells. We call the relative Cayley graphs of $G$ "the covers of $G$", because their geometric realizations are the 1-skeletons of the topological covers of $K$. Then $\mathrm{Cayley}(G, G_0)$ is a finite-sheeted cover (of the 1-skeleton of $K$) if and only if it has a finite number of vertices, which happens if and only if $G_0$ has finite index in $G$. However, the generating set $X$ of $G$ might be infinite, and then the finite-sheeted cover of $G$ is an infinite graph. To avoid possible conflicting terminology, we will not use the term "finite cover", and we say that a graph is finite if and only if it has finitely many vertices and edges.

*Definition* 2.2.   Let $E(\Gamma)$ denote the set of edges of a graph $\Gamma$. A labeling of $\Gamma$ by a set $X^*$ is a function $\mathrm{Lab} : E(\Gamma) \to X^*$ such that for any $e \in E(\Gamma)$, $\mathrm{Lab}(\bar{e}) = (\mathrm{Lab}(e))^{-1}$, where $\bar{e}$ denotes the inverse of the edge $e$.

A graph with a labeling function is called a labeled graph. Denote the set of all words in $X$ by $W(X)$, and denote the equality of two words by "$\equiv$". The label of a path $p = e_1 e_2 \cdots e_n$ in $\Gamma$, where $e_i \in E(\Gamma)$, is the word $\mathrm{Lab}(p) \equiv \mathrm{Lab}(e_1) \cdots \mathrm{Lab}(e_n) \in W(X)$.

*Definition* 2.3. Let $G$ be a group generated by a set $X$, let $\Gamma$ be a graph labeled with $X^*$, and let $p$ be a path in $\Gamma$. In this case, as usual, we identify the word $\mathrm{Lab}(p)$ with the corresponding element in $G$. Let $G_0$ be a subgroup of $G$. For any edge $(G_0 g, x)$ in $\mathrm{Cayley}(G, G_0)$ define $\mathrm{Lab}(G_0 g, x) = x$. Then for any path $p = (G_0 g, x_1)(G_0 g x_1, x_2) \cdots (G_0 g x_1 x_2 \cdots x_{n-1}, x_n)$ in $\mathrm{Cayley}(G, G_0)$, $\mathrm{Lab}(p) \equiv x_1 \cdots x_n \in W(X)$. Let $v_0$ be a vertex in $\Gamma$. Define $\mathrm{Lab}(\Gamma, v_0) = \{\mathrm{Lab}(p) | p \text{ is a loop in } \Gamma \text{ beginning at } v_0\}$.

*Remark* 2.4. It is easy to see that $\mathrm{Lab}(\Gamma, v_0)$ is a subgroup of $G$, and that $\mathrm{Lab}(\mathrm{Cayley}(G, G_0), G_0 \cdot 1) = G_0$.

*Definition* 2.5. We say that a connected subgraph $\Gamma$ of $\mathrm{Cayley}(G, S)$ represents $S$ and $g$, if $\Gamma$ contains $S \cdot 1$ and $S \cdot g$, and if $\mathrm{Lab}(\Gamma, S \cdot 1) = S$. We say that $\Gamma$ represents $S$, if $\Gamma$ contains $S \cdot 1$ and if $\mathrm{Lab}(\Gamma, S \cdot 1) = S$.

The following result from [Gi 1] shows a connection between the profinite topology and relative Cayley graphs.

THEOREM 2.6. *A finitely generated subgroup $S$ of $G$ is closed in $\mathrm{PT}(G)$ if and only if for any $g \notin S$ there exists a finite subgraph $\Gamma$ of $\mathrm{Cayley}(G, S)$ representing $S$ and $g$, which can be embedded in a cover of $G$ with finitely many vertices.*

In this paper we apply Theorem 2.6 to amalgamated free products of groups.

*Definition* 2.7. We denote the initial and the terminal vertices of $p$ by $\iota(p)$ and by $\tau(p)$ respectively, and the inverse of $p$ by $\bar{p}$.

*Definition* 2.8. Let $X^*$ and $Y^*$ be disjoint sets, and let $\Gamma$ be a graph labeled with $X^* \cup Y^*$. We say that a vertex $v$ in $\Gamma$ is bichromatic if there exist edges $e_1$ and $e_2$ in $\Gamma$ with $\iota(e_1) = \iota(e_2) = v, \mathrm{Lab}(e_1) \in X^*$ and $\mathrm{Lab}(e_2) \in Y^*$; otherwise we say that $v$ is monochromatic. We say that $\Gamma$ is monochromatic if the labels of all its edges are only in $X^*$ or only in $Y^*$. An $X^*$-component of $\Gamma$ is a maximal connected subgraph of $\Gamma$ labeled with $X^*$, which contains at least one edge. A $Y^*$-component of $\Gamma$ is a maximal connected subgraph of $\Gamma$ labeled with $Y^*$, which contains at least one edge.

*Definition* 2.9. Let $X$ be a generating set of a group $G$ and let $Y$ be a generating set of a group $H$, such that $X^* \cap Y^* = \emptyset$. Let $\phi$ be an isomorphism between the subgroups $G_0$ of $G$ and $H_0$ of $H$, and let $A = G \underset{G_0 = H_0}{*} H$ be the amalgamated free product of $G$ and $H$ defined by $\phi$. We say that a subgraph $\Gamma$ of a relative Cayley graph of $A$ is a precover of $A$, if each $X^*$-component of $\Gamma$ is a cover of $G$ and each $Y^*$-component of $\Gamma$ is a cover of $H$.

In order to show that a finitely generated subgroup $S$ of $A$ is closed in $\mathrm{PT}(A)$, for any $a \notin S$ we choose a finite subgraph $\Gamma$ of $\mathrm{Cayley}(A, S)$ representing $S$ and $a$, and try to embed it in a precover of $A$ with finitely many vertices (if $S$ is finitely generated, then a finite graph representing $S$ and $a$ can be easily constructed; cf. [Gi 2]). Then we try to embed such a precover in a cover of $A$ with finitely many vertices.

If $G$ and $H$ are LERF, then any monochromatic component of such $\Gamma$ can be embedded in a cover of $G$ or of $H$ with finitely many vertices; so we can embed $\Gamma$ in a graph $\Gamma'$ with finitely many vertices such that each monochromatic component of $\Gamma'$ is a cover of $G$ or of $H$. We would like to know when such $\Gamma'$ is a precover or a cover of $A$.

*Definition* 2.10.  Let $\Gamma$ be a graph labeled with a set $S^*$ and let $S_0 \subset S^*$. Following [G-T] we say that $\Gamma$ is $S_0$-saturated at a vertex $v$, if for any $s \in S_0$ there exists $e \in E(\Gamma)$ with $\iota(e) = v$ and $\mathrm{Lab}(e) = s$. We say that $\Gamma$ is $S_0$-saturated, if it is $S_0$-saturated at any $v \in V(\Gamma)$.

*Definition* 2.11.  Let $A = G \underset{G_0 = H_0}{*} H$ be as in Definition 2.9. We say that a graph $\Gamma$ labeled with $X^* \cup Y^*$ is $A$-compatible at a bichromatic vertex $v$, if for any pair of monochromatic paths of different colors $p$ and $q$ in $\Gamma$ such that $\iota(p) = v = \iota(q)$, if $\mathrm{Lab}(p) = \mathrm{Lab}(q) \in G_0$, then $\tau(q) = \tau(p)$. We say that $\Gamma$ is $A$-compatible, if it is $A$-compatible at all bichromatic vertices.

The following result from [Gi 2] gives a characterization of covers and precovers of $A$.

LEMMA 2.12.  *Let $\Gamma$ be a graph labeled with $X^* \cup Y^*$ such that each $X^*$-component of $\Gamma$ is a cover of $G$ and each $Y^*$-component of $\Gamma$ is a cover of $H$. Then $\Gamma$ is a precover of $A$ if and only if $\Gamma$ is $A$-compatible, and $\Gamma$ is a cover of $A$ if and only if, in addition, $\Gamma$ is $(X^* \cup Y^*)$-saturated.*

In the special case when the amalgamated free product is a double, i.e. the map $\phi$ in Definition 2.9 is the restriction of an isomorphism $\alpha$ from $G$ to $H$ (see Definition 1.1), the following result from [Gi 1] emphasizes the importance of precovers. We include the proof, as [Gi 1] is not easily available.

THEOREM 2.13 (the doubling theorem).  *Let $D$ be the double of a group $G$ along a subgroup $G_0$. Then any precover $\Gamma$ of $D$ with finitely many vertices can be embedded in a cover of $D$ with finitely many vertices.*

*Proof.* Define a new precover $\bar{\Gamma}$ of $D$ as follows. Let $\bar{\Gamma}$ be an abstract unlabeled graph isomorphic to $\Gamma$ and let $\beta : \bar{\Gamma} \to \Gamma$ be an isomorphism. For any edge $e$ of $\bar{\Gamma}$ define $\text{Lab}(e) = \alpha(\text{Lab}(\beta(e)))$ if $\text{Lab}(\beta(e)) \in X^*$, and $\text{Lab}(e) = \alpha^{-1}(\text{Lab}(\beta(e)))$ if $\text{Lab}(\beta(e)) \in Y^*$, where $\alpha, X$ and $Y$ are as in Definition 1.1. Then $\bar{\Gamma}$ is labeled with $X^* \cup Y^*$, and Lemma 2.12 implies that $\bar{\Gamma}$ is a precover of $D$. Indeed, as $\alpha$ and $\beta$ are isomorphisms, each monochromatic component of $\bar{\Gamma}$ is a cover of $G$ or of $H$. Let $v$ be a bichromatic vertex in $\bar{\Gamma}$, and let $p$ and $q$ be monochromatic paths of different colors in $\bar{\Gamma}$ such that $\iota(p) = v = \iota(q)$ and $\text{Lab}(p) = \text{Lab}(q) \in G_0$. Then $\beta(v)$ is a bichromatic vertex in $\Gamma$, and $\beta(p)$ and $\beta(q)$ are monochromatic paths of different colors in $\Gamma$ such that $\iota(\beta(p)) = \beta(v) = \iota(\beta(q))$ and $\text{Lab}(\beta(p)) = \text{Lab}(\beta(q)) \in G_0$. As $\Gamma$ is a precover, it is $D$-compatible at $\beta(v)$; hence $\tau(\beta(q)) = \tau(\beta(p))$, but then $\tau(q) = \tau(p)$ and therefore $\bar{\Gamma}$ is $D$-compatible.

Let $\Gamma'$ be a graph constructed from the disjoint union of $\Gamma$ and $\bar{\Gamma}$ by identifying every monochromatic vertex $v \in V(\bar{\Gamma})$ with $\beta(v) \in V(\Gamma)$. Then $\Gamma'$ has finitely many vertices and $\Gamma$ is embedded in $\Gamma'$. As $\Gamma$ and $\bar{\Gamma}$ are precovers, each monochromatic component of $\Gamma'$ is a cover of $G$ or of $H$. Let $v'$ be a bichromatic vertex in $\Gamma'$, and let $p'$ and $q'$ be monochromatic paths of different colors in $\Gamma'$ such that $\iota(p') = v' = \iota(q')$ and $\text{Lab}(p') = \text{Lab}(q') \in G_0$. If $v'$ has a preimage in $\Gamma$ which is bichromatic in $\Gamma$, then as each monochromatic component of $\Gamma$ is a cover of $G$ or of $H$, $p'$ and $q'$ have unique preimages in $\Gamma$. As $\Gamma$ is $D$-compatible at the preimage of $v'$, the preimages of $p'$ and $q'$ in $\Gamma$ have the same terminal vertex, but then $p'$ and $q'$ have the same terminal vertex in $\Gamma'$. The same argument shows that $\Gamma'$ is $D$-compatible at $v'$ if $v'$ has a preimage in $\bar{\Gamma}$ which is bichromatic in $\bar{\Gamma}$. If the preimage of $v'$ in $\Gamma$ is monochromatic, then $v'$ also has a monochromatic preimage in $\bar{\Gamma}$, so one path, say $p'$, has a unique preimage $p$ in $\Gamma$ and the other, $q'$, has a unique preimage $q$ in $\bar{\Gamma}$. Note that the path $\beta(q)$ belongs to $\Gamma$, $\text{Lab}(\beta(q)) = \text{Lab}(q) = \text{Lab}(q') = \text{Lab}(p') = \text{Lab}(p)$ and $\iota(\beta(q)) = \beta(\iota(q)) = \iota(p)$. Hence as $\Gamma$ is a precover, $\tau(\beta(q)) = \tau(p)$. As $\beta(\tau(q)) = \tau(\beta(q))$, the definition of $\Gamma'$ implies that $\tau(p') = \tau(q')$, so that $\Gamma'$ is $D$-compatible at $v'$. As $\Gamma'$ is $X^* \cup Y^*$-saturated, Lemma 2.12 implies that $\Gamma'$ is a cover of $D$.

Let $M_G$ and $M_H$ be topological manifolds of the same dimension, and let $M_{G_0}$ and $M_{H_0}$ be isomorphic boundary components of $M_G$ and $M_H$, respectively. Let $M_A$ be the manifold constructed from the disjoint union of $M_G$ and $M_H$ by identifying $M_{G_0}$ and $M_{H_0}$ via the fixed isomorphism. The concept of a precover can be restated in this category, and then the proof of the doubling theorem has an obvious geometrical interpretation. In fact, the concept of a precover and the doubling theorem can be restated for any pair $T_G$ and $T_H$ of topological spaces and their isomorphic subspaces $T_{G_0}$ and $T_{H_0}$.

Theorems 2.6 and 2.13 provide an important characterization of subgroups closed in the profinite topology on doubles.

COROLLARY 2.14. *A finitely generated subgroup $S$ is closed in $\mathrm{PT}(D)$ if and only if for any $d \notin S$ there exists a finite subgraph of $\mathrm{Cayley}(D, S)$, representing $S$ and $d$, which can be embedded in a precover of $D$ with finitely many vertices.*

Definition 2.15. A labeled graph is called well-labeled if for any $e_1$ and $e_2$ in $E(\Gamma)$ with $\iota(e_1) = \iota(e_2)$, if $\mathrm{Lab}(e_1) = \mathrm{Lab}(e_2)$, then $\tau(e_1) = \tau(e_2)$.

The following result from [Gi 2] will be used in the proof of Theorem 5.4.

LEMMA 2.16. *A graph $\Gamma$, well-labeled with the set $X^*$, can be embedded in a cover of $G$ if and only if any path $p$ in $\Gamma$ with $\mathrm{Lab}(p) = 1$ is a loop, i.e. $\iota(p) = \tau(p)$.*

In this paper we use the special case of the amalgamation of graphs ([Sta], [Gi 2]), which we call "grafting".

Definition 2.17. Let $G_0$ be a subgroup of $G$. Choose generating sets $X^*$ for $G$ and $X_1^*$ for $G_0$ such that $X_1^* \subset X^*$. Let $\Gamma$ be a graph well-labeled with $X^*$, and let $\beta_v$ be the $X_1^*$-component of the vertex $v$ in $\Gamma$. Let $\alpha$ be a graph well-labeled with $X_1^*$ such that $(\beta, v)$ embeds in $(\alpha, w)$. The graft of $(\alpha, v)$ on $(\Gamma, w)$ is constructed by taking the disjoint union of $\alpha$ and $\Gamma$, identifying the vertices $v$ and $w$, and then identifying two copies of $(\beta, v)$.

LEMMA 2.18. *The graft $\Delta$ of $\alpha$ on $\Gamma$ is well-labeled with $X^*$, and $\alpha$ and $\Gamma$ imbed in $\Delta$.*

*Proof.* Let $e_1$ and $e_2$ be edges in $\Delta$ with $\mathrm{Lab}(e_1) = \mathrm{Lab}(e_2)$ and $\iota(e_1) = \iota(e_2)$. If both $e_1$ and $e_2$ are in $\alpha$, then $e_1 = e_2$, because $\alpha$ is well-labeled with $X_1^*$. If both $e_1$ and $e_2$ are in $\Gamma$, then $e_1 = e_2$, because $\Gamma$ is well-labeled with $X^*$. If one edge, say $e_1$, is in $\alpha$, and another is in $\Gamma$, then $\mathrm{Lab}(e_1) = \mathrm{Lab}(e_2) \in X_1^*$ and $\iota(e_1) \in \Gamma \cap \alpha$. Hence $\iota(e_1) \in \beta$, but then, as $\beta$ is an $X_1^*$-component in $\Gamma$, $e_1 \in \beta \subset \Gamma$ and $e_2 \in \beta \subset \alpha$. Therefore by construction of $\Delta$, $e_1 = e_2$, so that $\Delta$ is well-labeled with $X^*$. By definition of grafting, we do not identify edges of $\Gamma$ with each other or edges of $\alpha$ with each other; hence $\Gamma$ and $\alpha$ are embedded in $\Delta$.

Note that, in general, graphs do not embed in their amalgams ([Sta], [Gi 2]).

## 3. Constructions of precovers

All the results in this section are valid for any amalgamated free product $A = G *_{G_0 = H_0} H$ (and not only for a double of $G$), and Lemma 3.1 holds for any groups $G$ and $H$ (they do not have to be LERF or negatively curved).

LEMMA 3.1. *Let $\Gamma$ be a graph with finitely many vertices which has the following properties.*

1) *All monochromatic components of $\Gamma$ are covers of $G$ or of $H$;*
2) *For any bichromatic vertex $v$ of $\Gamma$, $\mathrm{Lab}(\Gamma^{X^*}, v) \cap G_0 = \mathrm{Lab}(\Gamma^{Y^*}, v) \cap G_0$, where $\Gamma^{X^*}$ and $\Gamma^{Y^*}$ are, respectively, the $X^*$-component and the $Y^*$-component of $\Gamma$ containing $v$;*
3) *For any pair of bichromatic vertices in $\Gamma$ connected by a monochromatic path $p$ labeled by an element in $G_0$, there exists a pair $p'$ and $q'$ of monochromatic paths of different colors with the same endpoints as $p$, such that $\mathrm{Lab}(p') = \mathrm{Lab}(q') \in G_0$.*

*Then $\Gamma$ can be mapped onto a precover $\Pi$ of $A$ with finitely many vertices, by identifying certain pairs of monochromatic vertices of different color. This mapping restricts to an embedding on the union of all monochromatic components of the same color.*

*Proof.* If $\Gamma$ is $A$-compatible, then Lemma 2.12 implies that $\Gamma$ is a precover. Otherwise, there exists a bichromatic vertex $v$ in $\Gamma$ and monochromatic paths $p$ and $q$ of different colors in $\Gamma$ which begin at $v$, such that $\mathrm{Lab}(p) = \mathrm{Lab}(q) \in G_0$, but $\tau(p) \neq \tau(q)$. This might happen only if $\tau(p)$ and $\tau(q)$ are monochromatic vertices of different colors. Indeed, without loss of generality assume that $\tau(p)$ is bichromatic, then property 3 of $\Gamma$ implies that there exist a monochromatic path $p'$ of the same color as $p$ and a monochromatic path $q'$ of a different color, such that $p, p'$ and $q'$ have the same endpoints and $\mathrm{Lab}(p') = \mathrm{Lab}(q') \in G_0$. Then the path $p\bar{p}'$ is monochromatic and $\mathrm{Lab}(p\bar{p}') \in G_0$; hence property 2 of $\Gamma$ implies that there exists a closed monochromatic path $q''$ of the same color as $q$ with $\iota(q'') = \iota(p)$ and $\mathrm{Lab}(q'') = \mathrm{Lab}(p\bar{p}')$. Then $\mathrm{Lab}(\bar{q}q''q') = \mathrm{Lab}(\bar{q})\mathrm{Lab}(p\bar{p}')\mathrm{Lab}(q') = 1$; hence property 1 of $\Gamma$ implies that $\bar{q}q''q'$ is a closed loop, and so $q$ has the same endpoints as $q'$. Hence $\tau(p) = \tau(q)$, a contradiction.

Also for any vertex $u$ in $\Gamma$ there exists at most one vertex $w \neq u$ with the following property: there exists a pair of monochromatic paths $t$ and $s$ of different colors in $\Gamma$ such that $\tau(t) = u, \tau(s) = w, \iota(t) = \iota(s)$ and $\mathrm{Lab}(t) = \mathrm{Lab}(s) \in G_0$. Indeed, assume that there exists a vertex $w' \neq w$ and corresponding paths $t'$ and $s'$. If $\mathrm{Lab}(t) = \mathrm{Lab}(t')$, then property 1 of $\Gamma$ implies that $t'\bar{t}$ is a closed path. Then $\iota(s') = \iota(s)$ and $\mathrm{Lab}(s) = \mathrm{Lab}(s')$; so property 1 of $\Gamma$ implies that $s'\bar{s}$ is a closed path, hence $w = w'$.

If $\mathrm{Lab}(t) \neq \mathrm{Lab}(t')$, then $t'\bar{t}$ is a monochromatic path labeled with an element in $G_0$ which joins the initial vertices of $t'$ and $t$. Hence property 3 of $\Gamma$ implies that there exist monochromatic paths $t''$ and $s''$ of different colors in $\Gamma$ joining $\iota(t')$ to $\iota(t)$ such that $\mathrm{Lab}(t'') = \mathrm{Lab}(s'') \in G_0$, and such that $t''$ has the same color as $t$.

But then $t''t\bar{t}'$ is a monochromatic closed loop with $\mathrm{Lab}(t''t\bar{t}') \in G_0$; hence property 2 of $\Gamma$ implies that there exists a monochromatic loop $s_0$ of the same color as $s$, with the same initial vertex and the same label as $t''t\bar{t}'$. But then $\mathrm{Lab}(\bar{s}\bar{s}''s_0s') = \mathrm{Lab}(\bar{t})\mathrm{Lab}(\bar{t}'')\mathrm{Lab}(t''t\bar{t}')\mathrm{Lab}(t') = 1$; so property 1 of $\Gamma$ implies that $\bar{s}\bar{s}''s_0s'$ is a closed path, and thus $w = w'$.

We construct the mapping of $\Gamma$ onto a precover as follows. For any pair of monochromatic paths of different colors in $\Gamma$ which have the same label and the same initial vertex, but distinct terminal vertices, we identify their terminal vertices. As $\Gamma$ has finitely many vertices, after repeating this procedure a finite number of times, we obtain an $A$-compatible graph $\Pi$. The monochromatic components of $\Gamma$ coincide with the monochromatic components of $\Pi$, because the above discussion shows that we identify any monochromatic vertex in $\Gamma$ with at most one monochromatic vertex of different color, and the identifications do not involve bichromatic vertices. Hence property 1 of $\Gamma$ and Lemma 2.12 imply that $\Pi$ is a precover of $A$.

Let $\phi$ be as in Definition 2.9. To make the rest of the exposition easier to follow, we assume that the generating set $X_1$ of $G_0$ is a subset of $X$ and its image $Y_1 = \phi(X_1)$, which is a generating set of $H_0 = \phi(G_0)$, is a subset of $Y$.

*Remark* 3.2. Let $S$ be a finitely generated subgroup of $A$, let $a$ be an element in $A$, but not in $S$, and let $\Gamma'$ be a finite subgraph of $\mathrm{Cayley}(A, S)$ representing $S$ and $a$. Let $x_1$ be an element in $X_1$, and let $y_1 = \phi(x_1)$. For any vertex $v$ in $\Gamma'$, let $e_{v,x}$ and $e_{v,y}$ be edges in $\mathrm{Cayley}(A, S)$ which begin at $v$ and are labeled with $x_1$ and $y_1$, respectively. Define $\Gamma''$ to be the union of $\Gamma'$ and all the edges $e_{v,x}$ and $e_{v,y}$. Note that the edges $e_{v,x}$ and $e_{v,y}$ have the same terminal vertex; hence $\Gamma''$ is a finite subgraph of $\mathrm{Cayley}(D, S)$ representing $S$ and $a$, and all the vertices of $\Gamma''$ are bichromatic in $\Gamma''$. If we can embed $\Gamma''$ in a graph $\Gamma$ which has properties 1–3 of Lemma 3.1, then we can map $\Gamma$ onto a precover $\Pi$, as in Lemma 3.1 and, as all the vertices of $\Gamma''$ are bichromatic in $\Gamma$, this map of $\Gamma$ restricts to an embedding on $\Gamma''$. However, examples discussed in Section 2 show that such embeddings do not exist for arbitrary groups $S$ and $A$; otherwise any double of a LERF group would be LERF. The following result shows that under certain assumptions on $S$, we can almost achieve this goal.

LEMMA 3.3. *Let $S$ be a finitely generated subgroup of $A = G \underset{G_0 = H_0}{*} H$, such that the intersection of $S$ with any conjugate of $G_0$ in $A$ is finitely generated, and let $\phi$ be as in Definition* 2.9. *Then any finite subgraph $\Gamma_0$ of* Cayley$(A, S)$ *representing $S$, is contained in a finite subgraph $\Gamma_1$ of* Cayley$(A, S)$ *with the following properties*:

4) *For any bichromatic vertex $w$ of $\Gamma_1$,* Lab$(\Gamma_1^{X_1^*}, w) \cap G_0 =$ Lab$(\Gamma_1^{Y_1^*}, w) \cap G_0 =$ Lab(Cayley$(A,S), w) \cap G_0$, *where $\Gamma_1^{X_1^*}$ and $\Gamma_1^{Y_1^*}$ are, respectively, the $X_1^*$-component and the $Y_1^*$-component of $\Gamma_1$ containing $w$.*

5) *If two distinct bichromatic vertices in $\Gamma_1$ are connected by a path in* Cayley$(A, S)$ *labeled with an element of $G_0$, then they are connected by a pair of monochromatic paths $p'$ and $q'$ in $\Gamma_1$ labeled with $X_1^*$ and $Y_1^*$, respectively, such that $\phi($Lab$(p')) \equiv$ Lab$(q') \in G_0$.*

*Therefore $\Gamma_1$ has properties* 2 *and* 3 *of Lemma* 3.1.

*Proof.* Let $\Gamma_0$ be any finite subgraph of Cayley$(A, S)$ representing $S$. If $\Gamma_0$ already has properties 4 and 5, take $\Gamma_1 = \Gamma_0$. Otherwise, let $W$ be the set of all bichromatic vertices of $\Gamma_0$. For each pair of distinct vertices in $W$ which are connected by a path in Cayley$(A, S)$ labeled with an element of $G_0$, choose a pair of paths $p_0$ and $q_0$ in Cayley$(A, S)$ connecting these vertices, labeled with $X_1^*$ and $Y_1^*$ respectively, such that Lab$(q_0) \equiv \phi($Lab$(p_0))$.

For any $w \in W$, the group Lab(Cayley$(A, S), w$) is a conjugate of Lab(Cayley$(A, S), S \cdot 1) = S$; hence the subgroup Lab(Cayley$(A, S), w) \cap G_0$ is finitely generated (because it is a conjugate of the intersection of a conjugate of $G_0$ with $S$). Therefore we can choose a finite number of loops $p_{w,i}$ and $q_{w,i}$ in Cayley$(A, S)$ labeled with $X_1^*$ and $Y_1^*$ respectively, which begin at $w$ such that Lab$(q_{w,i}) \equiv \phi($Lab$(p_{w,i}))$, and such that the set $\{$Lab$(p_{w,i})\}$ (hence the set $\{$Lab$(q_{w,i})\}$) generates the subgroup Lab(Cayley$(A, S), w) \cap G_0$.

Let $\Gamma_1$ be the union of $\Gamma_0$ and all the paths $p_0, q_0, p_{w,i}$ and $q_{w,i}$. Then $\Gamma_1$ is a finite graph, and we will show that it has properties 4 and 5. By construction, $\Gamma_1$ has the required properties for all vertices in $W$. However, the set of bichromatic vertices of $\Gamma_1$ is bigger than $W$, as all the new vertices which were added to $\Gamma_0$ to construct $\Gamma_1$ are bichromatic in $\Gamma_1$. Hence for any bichromatic vertex $u \notin W$ in $\Gamma_1$ there exists a vertex $w \in W$ and paths $c$ and $d$ in $\Gamma_1$, labeled with $X_1^*$ and $Y_1^*$, respectively, joining $u$ to $w$ such that Lab$(d) \equiv \phi($Lab$(c)) \in G_0$.

Consider a path $p$ in Cayley$(A, S)$ joining distinct bichromatic vertices $v_1$ and $v_2$ of $\Gamma_1$, such that Lab$(p) \in G_0$. As was mentioned above, there exist paths $c_i$ and $d_i$ in $\Gamma_1$, labeled with $X_1^*$ and $Y_1^*$, respectively, such that Lab$(d_i) \equiv \phi($Lab$(c_i)) \in G_0$, and $c_i$ and $d_i$ join $v_i$ to some $w_i \in W, i = 1, 2$. Then $\bar{c}_1 p c_2$ is a path in Cayley$(A, S)$ labeled with an element in $G_0$ joining $w_1$ to $w_2$. As $\Gamma_1$ has property 5 for all vertices in $W$, there exists a pair of paths $c_0$ and $d_0$ labeled with $X_1^*$ and $Y_1^*$, respectively, in $\Gamma_1$ joining $w_1$ to $w_2$ such that $\phi(\mathrm{Lab}(c_0)) \equiv \mathrm{Lab}(d_0)$. Then $c_1 c_0 \bar{c}_2$ and $d_1 d_0 \bar{d}_2$ are paths in $\Gamma_1$ labeled by $X_1^*$ and $Y_1^*$, respectively, joining $v_1$ to $v_2$, and $\phi(\mathrm{Lab}(\bar{c}_1 c_0 c_2)) \equiv \mathrm{Lab}(d_1 d_0 \bar{d}_2)$; hence $\Gamma_1$ has property 5 for any pair of bichromatic vertices.

Consider a bichromatic vertex $v$ in $\Gamma_1$, and let $c$ and $d$ be paths labeled with $X_1^*$ and $Y_1^*$, respectively, in $\Gamma_1$ with $\mathrm{Lab}(d) \equiv \phi(\mathrm{Lab}(c)) \in G_0$, which join $v$ to some $w \in W$. Then $v$ and $w$ belong to the same $X_1^*$-component of $\Gamma_1$, say $\Gamma^{X_1^*}$, and to the same $Y_1^*$-component of $\Gamma_1$, say $\Gamma^{Y_1^*}$. Hence $\mathrm{Lab}(\Gamma^{X_1^*}, v) \cap G_0 = (\mathrm{Lab}(c)\mathrm{Lab}(\Gamma^{X_1^*}, w)\mathrm{Lab}(\bar{c})) \cap G_0 = \mathrm{Lab}(c)(\mathrm{Lab}(\Gamma^{X_1^*}, w) \cap G_0)\mathrm{Lab}(\bar{c})$. But $\Gamma_1$ has property 4 for any $w \in W$; hence $\mathrm{Lab}(\Gamma^{X_1^*}, w) \cap G_0 = \mathrm{Lab}(\Gamma^{Y_1^*}, w) \cap G_0 = \mathrm{Lab}(\mathrm{Cayley}(A,S), w) \cap G_0$. Therefore,

$$\begin{aligned}\mathrm{Lab}(\Gamma^{X_1^*}, v) \cap G_0 &= \mathrm{Lab}(c)(\mathrm{Lab}(\mathrm{Cayley}(A,S), w) \cap G_0)\mathrm{Lab}(\bar{c}) \\ &= \mathrm{Lab}(\mathrm{Cayley}(A,S), v) \cap G_0 \\ &= \mathrm{Lab}(d)(\mathrm{Lab}(\Gamma^{Y_1^*}, w) \cap G_0)\mathrm{Lab}(\bar{d}) = \mathrm{Lab}(\Gamma^{Y_1^*}, v) \cap G_0.\end{aligned}$$

So $\Gamma_1$ has property 4 for any bichromatic vertex.

Lemma 3.3 shows that under certain assumptions we can embed the graph $\Gamma''$ described in Remark 3.2 in a graph $\Gamma_1$ which has properties 2 and 3 of Lemma 3.1. However, our goal is to embed $\Gamma''$ in a graph $\Gamma$ which has properties 1–3 of Lemma 3.1. Unlike $\Gamma_1$, generically, such $\Gamma$ cannot be a subgraph of $\mathrm{Cayley}(A,S)$. It will be constructed using "grafting". The following lemma shows that a construction of a graph which has property 1 of Lemma 3.1 can be reduced to a construction of a graph which has two additional properties, which are easier to verify.

*Definition* 3.4.   Let $G_0$ be a subgroup of $G$. We say that a subgraph $\Gamma$ of a cover of $G$ is $G_0$-complete at a vertex $v$, if for any $g \in G_0$ there exists a path $p_g$ in $\Gamma$ beginning at $v$ with $\mathrm{Lab}(p_g) = g$.

LEMMA 3.5 (the grafting lemma).   *Let $G$ and $H$ be LERF groups. Let $\Sigma_0$ be a finite graph with the following properties.*

6) *All monochromatic components of $\Sigma_0$ are subgraphs of covers of $G$ or of $H$;*

7) *For any bichromatic vertex $w$ of $\Sigma_0$, the monochromatic components $\Sigma_0^{X^*}$ and $\Sigma_0^{Y^*}$ of $\Sigma_0$ containing $w$ are, respectively, $G_0$-complete and $H_0$-complete at $w$.*

*If $\Sigma_0$ has properties 2 and 3 of Lemma 3.1, then it can be embedded in a graph $\Sigma$ which has properties 1–3 of Lemma 3.1.*

*Proof.* As each monochromatic component of $\Sigma_0$ is a finite subgraph of a cover of $G$ or of $H$, and as $G$ is LERF, each monochromatic component of $\Sigma_0$ can be embedded in a cover of $G$ or of $H$ with finitely many vertices. Let $\Sigma$ be the graph constructed from the disjoint union of $\Sigma_0$ and all these covers by identifying each monochromatic component of $\Sigma_0$ with its image in the corresponding cover. (Here we use "grafting"). Then, by construction, $\Sigma$ has property 1 of Lemma 3.1.

Let $w$ be a bichromatic vertex in $\Sigma$, let $\Sigma^{X^*}$ be the $X^*$-component of $\Sigma$ containing $w$, and let $l$ be a loop in $\Sigma^{X^*}$ which begins at $w$ such that $\text{Lab}(l) \in G_0$. As $\Sigma$ and $\Sigma_0$ have the same sets of bichromatic vertices, $w$ is bichromatic in $\Sigma_0$. As the $X^*$-component $\Sigma_0^{X^*}$ of $\Sigma_0$ containing $w$ is $G_0$-complete at $w$, there exists a path $l'$ in $\Sigma_0^{X^*}$ which begins at $w$ with $\text{Lab}(l') = \text{Lab}(l)$. As $\Sigma_0^{X^*}$ is embedded in $\Sigma^{X^*}$, and as $\Sigma^{X^*}$ is a cover of $G$, the paths $l$ and $l'$ have the same terminal vertex (because they have the same initial vertex and the same label). Therefore $l'$ is a loop in $\Sigma_0^{X^*}$. As $\Sigma_0$ has property 2 of Lemma 3.1, the $Y^*$-component $\Sigma_0^{Y^*}$ of $\Sigma_0$ containing $w$ contains a loop $l''$ which begins at $w$ with $\text{Lab}(l'') = \text{Lab}(l')$; hence $\text{Lab}(\Sigma^{X^*}, w) \cap G_0$ is contained in $\text{Lab}(\Sigma^{Y^*}, w) \cap G_0$. Similarly, $\text{Lab}(\Sigma^{Y^*}, w) \cap G_0$ is contained in $\text{Lab}(\Sigma^{X^*}, w) \cap G_0$; therefore $\Sigma$ has property 2 of Lemma 3.1.

Consider bichromatic vertices $w_1$ and $w_2$ in $\Sigma$ connected by a monochromatic path $p$ with $\text{Lab}(p) \in G_0$. As $\Sigma$ and $\Sigma_0$ have the same sets of bichromatic vertices, $w_1$ and $w_2$ are bichromatic in $\Sigma_0$. As each monochromatic component of $\Sigma_0$ is $G_0$-complete at $w_1$, there exists a monochromatic path $p_1$ in $\Sigma_0$ beginning at $w_1$, which has the same color and the same label as $p$. As each monochromatic component of $\Sigma$ is a cover, $\tau(p) = \tau(p_1) = w_2$. As $\Sigma_0$ has property 3 of Lemma 3.1, there exist monochromatic paths $p_0$ and $q_0$ in $\Sigma_0$ of different colors connecting $w_1$ and $w_2$ such that $\text{Lab}(p_0) = \text{Lab}(q_0) \in G_0$. As $p_0$ and $q_0$ lie in $\Sigma$, it follows that $\Sigma$ also has property 3 of Lemma 3.1, as required.

Our next goal is to construct a graph which has property 7 of Lemma 3.5. The following lemma shows that we can easily do it in a very special case. The general case is considered in Theorem 4.4.

LEMMA 3.6. *If $G_0$ is finitely generated, then the graph $\Gamma_1$, constructed in the proof of Lemma 3.3 is contained in a finite subgraph $\Gamma_2$ of $\text{Cayley}(A, S)$ which has properties 4 and 5 of Lemma 3.3 and, in addition, is $(X_1^* \cup Y_1^*)$-saturated at any bichromatic vertex $u$ such that $\text{Lab}(\text{Cayley}(A, S), u) \cap G_0$ has finite index in $G_0$.*

*Hence $\Gamma_2$ has property 7 of Lemma 3.5 for any bichromatic vertex $u$ such that $\text{Lab}(\text{Cayley}(A, S), u) \cap G_0$ has finite index in $G_0$.*

*Proof.* Note that the definition of the sets $X_1$ and $Y_1$ implies that the $X_1^*$-component and $Y_1^*$-component of any vertex $u$ in Cayley$(A, S)$ are isomorphic covers of $G_0$, hence they are $G_0$-complete at any vertex. Also the sets of vertices of these components coincide, so that the union of the $X_1^*$-component and the $Y_1^*$-component of any vertex in Cayley$(A, S)$ consists entirely of vertices bichromatic in this union.

Let $U = \{u_1, \cdots, u_k\}$ be the set of all bichromatic vertices of $\Gamma_1$, such that Lab(Cayley$(A, S), u_i) \cap G_0$ is of finite index in $G_0$. Then the $X_1^*$-component and $Y_1^*$-component of any $u_i \in U$ in Cayley$(A, S)$ are finite. Define $\Gamma_2$ to be the union of $\Gamma_1$ and the $X_1^*$-components and the $Y_1^*$-components of all $u_i \in U$ in Cayley$(A, S)$. Then $\Gamma_2$ is a finite graph, which has property 7 of Lemma 3.5 at any vertex $u_i$. By construction of $\Gamma_2$, if $u$ is a bichromatic vertex in $\Gamma_2$ such that Lab(Cayley$(A, S), u) \cap G_0$ has finite index in $G_0$, then $u$ belongs to the $X_1^*$-component (and to the $Y_1^*$-component) of some $u_i \in U$; hence $\Gamma_2$ has property 7 of Lemma 3.5 at any such $u$. It is easy to see that $\Gamma_2$ has properties 4 and 5 of Lemma 3.3, because $\Gamma_1$ has them.

COROLLARY 3.7. *Let $G_0$ be finitely generated. A special case of Lemma 3.6 with $H = H_0$ states that for any finitely generated subgroup $S$ of $G$, such that the intersection of $S$ with any conjugate of $G_0$ in $G$ is finitely generated, and for any finite subgraph $\Gamma$ of Cayley$(G, S)$ representing $S$, there exists a finite subgraph $\Gamma'$ of Cayley$(G, S)$ containing $\Gamma$ with the following properties:*

4′) *For any vertex $w$ of $\Gamma'$ and for any $g \in $ Lab(Cayley$(G, S), w) \cap G_0$ there exists a loop $l_g$ in the $X_1^*$-component of $w$ in $\Gamma'$ which begins at $w$, such that Lab$(l_g) = g$.*

5′) *Any two vertices in $\Gamma'$ are joined by a path in Cayley$(G, S)$ labeled by an element in $G_0$ if and only if they belong to the same $X_1^*$-component in $\Gamma'$. (The "if" direction always holds.)*

7′) *$\Gamma'$ is $X_1^*$-saturated (hence, it is $G_0$-complete) at any vertex $v$ such that Lab(Cayley$(G, S), v) \cap G_0$ has finite index in $G_0$.*

## 4. Strongly separable subgroups

We will use the following fact proved in [Gi 2].

LEMMA 4.1. *If $(\Delta, u)$ is a subgraph of a cover of a group $G$, then $(\Delta, u)$ can be isomorphically embedded in the relative Cayley graph (Cayley$(G, $ Lab$(\Delta, u)),$ Lab$(\Delta, u) \cdot 1)$. To avoid awkward notation, we denote this relative Cayley graph by $(\tilde{\Delta}, u)$.*

The notation (Cayley$(G, $ Lab$(\Delta, u)),$ Lab$(\Delta, u) \cdot 1) = (\tilde{\Delta}, u)$ will be used through the rest of the paper.

Let $G_0$ be a finitely generated subgroup of a group $G$. We choose generating sets $X_1$ of $G_0$ and $X$ of $G$, such that $X_1$ is a finite subset of $X$.

*Definition* 4.2. We use the notation of Lemma 4.1. We say that a finitely generated subgroup $G_0$ of a group $G$ is strongly separable in $G$, if for any finitely generated subgroups $S_1$ and $S_2$ of $G$ such that the subgroups $S_1 \cap G_0$ and $S_2 \cap G_0$ are equal and have infinite index in $G_0$, and for any finite subgraphs $(\Gamma_1, v_1)$ and $(\Gamma_2, v_2)$ of covers of $G$ with $\text{Lab}(\Gamma_1, v_1) = S_1$ and $\text{Lab}(\Gamma_2, v_2) = S_2$, there exist finite subgraphs $\Gamma'_i, (i = 1, 2)$ of $\tilde{\Gamma}_i$ which contain $\Gamma_i$ and have properties $4'$, $5'$ and $7'$ of Corollary 3.7 in $\tilde{\Gamma}_i$, and there exist embeddings $(\Gamma'_i, v_i) \to (\Delta_i, v_i)$ with the following properties.

a) $\Delta_i$ is a finite subgraph of a cover of $G$, such that the $X_1^*$-component of $v_i$ in $\Delta_i$ is a cover of $G_0$; hence $\Delta_i$ is $G_0$-complete at $v_i$. If an edge $e$ of $\Delta_i$ does not belong to $\Gamma'_i$, then $e$ belongs to the $X_1^*$-component of $v_i$ in $\Delta_i, i = 1, 2$.

b) $\text{Lab}(\Delta_1, v_1) \cap G_0 = \text{Lab}(\Delta_2, v_2) \cap G_0$.

c) If $u_i \in V(\Gamma'_i)$ does not belong to the $X_1^*$-component of $v_i$ in $(\tilde{\Gamma}_i, v_i)$, then $\text{Lab}(\tilde{\Gamma}_i, u_i) \cap G_0 = \text{Lab}(\tilde{\Delta}_i, u_i) \cap G_0, i = 1, 2$.

d) A pair of vertices in the image of $\Gamma'_i$ belongs to the same $X_1^*$-component in $\tilde{\Delta}_i$ if and only if they belong to the same $X_1^*$-component in $\tilde{\Gamma}_i$.

Let $S$ be a subgroup of a group $G$. Note that the vertices $Sg_1$ and $Sg_2$ of $\text{Cayley}(G, S)$ belong to the same $X_1^*$-component in $\text{Cayley}(G, S)$ if and only if $g_1 \in Sg_2G_0$; hence Definition 4.2 can be equivalently restated in pure group-theoretical language, but such a change of language would greatly complicate the proof of Theorem 5.4. Definition 4.2 looks very complicated, but it defines a nontrivial class of objects. Corollary 3.7 implies that a finite subgroup is strongly separable in any group, and Theorem 5.4 shows that certain groups have rich families of infinite strongly separable subgroups.

The following lemma demonstrates why the strong separability is useful in constructing embeddings of graphs.

LEMMA 4.3. *Let $D$ be a double of a LERF group $G$ along a finitely generated subgroup $G_0$, which is strongly separable in $G$. Let $\Gamma$ be a finite graph which has properties 2 and 3 of Lemma 3.1, and property 6 of Lemma 3.5. Assume that each monochromatic component $\Gamma_i$ of $\Gamma$ has properties $4'$, $5'$ and $7'$ of Corollary 3.7 in $\tilde{\Gamma}_i$. If $\Gamma$ does not have property 7 of Lemma 3.5, then $\Gamma$ can be embedded in a finite graph $\Omega$ which has the same set of bichromatic vertices and the same properties as $\Gamma$, and in addition, the set of all bichromatic vertices of $\Omega$, where $\Omega$ does not have property 7, is strictly smaller then the corresponding set in $\Gamma$.*

*Proof.* Let $u$ be a bichromatic vertex in $\Gamma$ and let $\Gamma^{X^*}$ and $\Gamma^{Y^*}$ be, respectively, the $X^*$-component and the $Y^*$-component of $\Gamma$ containing $u$. If $\Gamma^{X^*}$ is $G_0$-complete at $u$ then, as $\Gamma^{X^*}$ is a finite graph, it follows that $\text{Lab}(\Gamma^{X^*}, u) \cap G_0$ is of finite index in $G_0$. Property 2 of Lemma 3.1 states that $\text{Lab}(\Gamma^{X^*}, u) \cap G_0 = \text{Lab}(\Gamma^{Y^*}, u) \cap G_0$; hence $\text{Lab}(\Gamma^{Y^*}, u) \cap G_0$ is of finite index in $G_0$. As $\Gamma^{Y^*}$ has property $4'$ of Corollary 3.7, it follows that $\text{Lab}(\tilde\Gamma^{Y^*}, u) \cap G_0$ is of finite index in $G_0$. As $\Gamma^{Y^*}$ has property $7'$ of Corollary 3.7, it follows that $\Gamma^{Y^*}$ is $G_0$-complete at $u$. Therefore, if $\Gamma$ does not have property 7 of Lemma 3.5 at $u$, then both $\Gamma^{X^*}$ and $\Gamma^{Y^*}$ are not $G_0$-complete at $u$, and $\text{Lab}(\Gamma^{X^*}, u) \cap G_0 = \text{Lab}(\Gamma^{Y^*}, u) \cap G_0$ is of infinite index in $G_0$.

Then the strong separability of $G_0$ in $G$ implies the existence of a finite subgraph $\Gamma'^{X^*}$ of $\tilde\Gamma^{X^*}$ which contains $\Gamma^{X^*}$ and has properties $4'$, $5'$ and $7'$ of Corollary 4.3 in $\tilde\Gamma^{X^*}$, a finite subgraph $\Gamma'^{Y^*}$ of $\tilde\Gamma^{Y^*}$ which contains $\Gamma^{Y^*}$ and has properties $4'$, $5'$ and $7'$ in $\tilde\Gamma^{Y^*}$, and the embeddings $(\Gamma'^{X^*}, u) \to (\Delta^{X^*}, u)$ and $(\Gamma'^{Y^*}, u) \to (\Delta^{Y^*}, u)$, which have properties a–d of Definition 4.2.

Let $\Omega$ be the graph constructed from the disjoint union of $\Gamma, \Delta^{X^*}$ and $\Delta^{Y^*}$, by identifying $\Gamma'^{X^*}$ with its image in $\Delta^{X^*}$, and $\Gamma'^{Y^*}$ with its image in $\Delta^{Y^*}$. (Here we use "grafting" again.) As $\Delta^{X^*}$ and $\Delta^{Y^*}$ are finite graphs, so is $\Omega$. Lemma 2.18 states that the inclusion of $\Gamma$ into $\Omega$ is an embedding. By construction, the monochromatic components of $\Omega$ containing $u$ are $\Delta^{X^*}$ and $\Delta^{Y^*}$, and the remaining monochromatic components of $\Omega$ are isomorphic to the corresponding monochromatic components of $\Gamma$; hence $\Omega$ has property 6 of Lemma 3.5. As the $X_1^*$ component of $u$ in $\Delta^{X^*}$ and the $Y_1^*$ component of $u$ in $\Delta^{Y^*}$ are covers of $G_0$, the same is true in $\Omega$; hence both monochromatic components of $\Omega$ which contain $u$ are $G_0$-complete at $u$. As $\Delta^{X^*}$ and $\Delta^{Y^*}$ are monochromatic, $\Gamma$ and $\Omega$ have the same set of bichromatic vertices.

$\Omega$ has property 3 of Lemma 3.1. Indeed, consider a monochromatic path $p$ in $\Omega$ with $\text{Lab}(p) \in G_0$ joining bichromatic vertices $v$ and $w$. If $p$ belongs to $\Gamma$, the result follows because $\Gamma$ has property 3 of Lemma 3.1. So we may assume, without loss of generality, that $p$ belongs to $\Delta^{X^*}$. The definition of $X_1^*$ implies that $v$ and $w$ belong to the same $X_1^*$-component in $\tilde\Delta^{X^*}$. Then property d of Definition 4.1 implies that $v$ and $w$ belong to the same $X_1^*$-component in $\tilde\Gamma^{X^*}$, and property $5'$ of $\Gamma^{X^*}$ implies that $v$ and $w$ belong to the same $X_1^*$-component in $\Gamma^{X^1}$. But then property 3 of $\Gamma$ implies that $\Gamma$ contains paths $p'$ and $q'$ labeled with $X^*$ and $Y^*$, respectively, with $\text{Lab}(p') = \text{Lab}(q') \in G_0$. As $\Gamma$ is embedded in $\Omega$, the paths $p'$ and $q'$ belong to $\Omega$, Hence $\Omega$ has property 3 of Lemma 3.1.

Any monochromatic component $\Omega_i$ of $\Omega$ has property $5'$ of Corollary 3.7 in $\tilde\Omega_i$. Indeed, without loss of generality, it is enough to show that $\Delta^{X^*}$ has property $5'$ in $\tilde\Delta^{X^*}$. Let $v$ and $w$ be a pair of vertices in $\Delta^{X^*}$ which belong to the same $X_1^*$-component in $\tilde\Delta^{X^*}$. Then property d of Definition 4.2 implies that they belong to the same $X_1^*$-component in $\tilde\Gamma^{X^*}$, and property $5'$ of $\Gamma^{X^*}$ implies that they belong to the same $X_1^*$-component in $\Gamma^{X^*}$. As $\Gamma^{X^*}$ is embedded in $\Delta^{X^*}$, it follows that $\Delta^{X^*}$ has property $5'$.

Any monochromatic component $\Omega_i$ of $\Omega$ has property $4'$ of Corollary 3.7 in $\tilde{\Omega}_i$. Indeed, without loss of generality, it is enough to show that $\Delta^{X^*}$ has property $4'$ in $\tilde{\Delta}^{X^*}$. Let $\Delta^{X_1^*}$ be the $X_1^*$-component of $u$ in $\Delta^{X^*}$, and let $w$ be a vertex in $\Delta^{X^*}$. Consider two cases.

i) If $w$ belongs to $\Delta^{X_1^*}$, then as $\Delta^{X_1^*}$ is a cover of $G_0$, it is $G_0$-complete at $w$. So for any $g \in \text{Lab}(\tilde{\Delta}^{X^*}, w) \cap G_0$ there exists a path $l_g$ in $\Delta^{X_1^*}$ which begins at $w$ such that $\text{Lab}(l_g) = g$. But as $\tilde{\Delta}^{X^*}$ is a cover of $G$, the path $l_g$ should be a loop; hence property $4'$ holds at $w$.

ii) If $w$ does not belong to $\Delta^{X_1^*}$, then property a of Definition 4.2 implies that the $X_1^*$-components of $w$ in $\Delta^{X^*}$ and in $\Gamma^{'X^*}$ coincide. Property c of Definition 4.2 states that $\text{Lab}(\tilde{\Delta}^{X^*}, w) \cap G_0 = \text{Lab}(\tilde{\Gamma}^{X^*}, w) \cap G_0$. But then, as $\Gamma^{'X^*}$ has property $4'$ at $w$, so does $\Delta^{X^*}$.

Any monochromatic component $\Omega_i$ of $\Omega$ has property $7'$ of Corollary 3.7 in $\tilde{\Omega}_i$. Indeed, without loss of generality, it is enough to show that $\Delta^{X^*}$ has property $7'$. So let $w$ be a vertex in $\Delta^{X^*}$ such that $\text{Lab}(\tilde{\Delta}^{X^*}, w) \cap G_0$ has finite index in $G_0$. If $w$ belongs to $\Delta^{X_1^*}$, then as $\Delta^{X_1^*}$ is a cover of $G_0$, it is $X_1^*$-saturated at $w$.

If $w$ does not belong to $\Delta^{X_1^*}$, then property c of Definition 4.2 implies that $\text{Lab}(\tilde{\Gamma}^{X^*}, w) \cap G_0$ has finite index in $G_0$. But then, as $\Gamma^{'X^*}$ has property $7'$, it is $X_1^*$-saturated at $w$, so is $\Delta^{X^*}$.

$\Omega$ has property 2 of Lemma 3.1. Indeed, consider a bichromatic vertex $v$ in $\Omega$. If $v$ does not belong to either $\Delta^{X^*}$ or $\Delta^{Y^*}$, then $\Omega$ has property 2 of Lemma 3.1 at $v$, because all the monochromatic components of $\Gamma$ and $\Omega$, except for $\Delta^{X^*}$ and $\Delta^{Y^*}$, coincide, and $\Gamma$ has property 2 at any bichromatic vertex by assumption. Now assume, without loss of generality, that $v$ belongs to $\Delta^{X^*}$. Property b of Definition 4.2 implies that $\Omega$ has property 2 of Lemma 3.1 at $u$.

Consider two cases.

i) If $v$ belongs to $\Delta^{X_1^*}$, then there exists a path labeled with $X_1^*$ joining $v$ and $u$ in $\Delta^{X^*}$. As $\Omega$ has property 3 of Lemma 3.1, it follows that $\Omega$ contains paths $p'$ and $q'$ labeled with $X^*$ and $Y^*$, respectively, with $\text{Lab}(p') = \text{Lab}(q') \in G_0$, which join $u$ to $v$. But then conjugation by $\text{Lab}(p')$ and by $\text{Lab}(q')$ shows that $\Omega$ has property 2 at $v$, because it has property 2 at $u$.

ii) If $v$ does not belong to $\Delta^{X_1^*}$ then, as was shown above, $v$ does not belong to $\Delta^{Y_1^*}$. Then property a of Definition 4.2 implies that the $X_1^*$-components of $w$ in $\Delta^{X^*}$ and in $\Gamma^{'X^*}$ coincide, and the $Y_1^*$-components of $w$ in $\Delta^{Y^*}$ and in $\Gamma'^{Y^*}$ coincide. As

$$(\text{Lab}(\Gamma^{X_1^*}, v) \cap G_0) < (\text{Lab}(\Gamma'^{X_1^*}, v) \cap G_0) < (\text{Lab}(\Gamma'^{X^*}, v) \cap G_0)$$
$$< (\text{Lab}(\tilde{\Gamma}^{X^*}, v) \cap G_0),$$

property 4' of $\Gamma^{X^*}$ implies that

$$\text{Lab}(\Gamma'^{X_1^*}, v) \cap G_0 = \text{Lab}(\Gamma^{X^*}, v) \cap G_0.$$

Similarly, property 4' of $\Delta^{X^*}$ implies that

$$\text{Lab}(\Delta^{X_1^*}, v) \cap G_0 = \text{Lab}(\Delta^{X^*}, v) \cap G_0.$$

Hence $\text{Lab}(\Gamma^{X^*}, v) \cap G_0 = \text{Lab}(\Delta^{X^*}, v) \cap G_0$. The corresponding equality holds for the $Y^*$-components; hence $\Omega$ has property 2 of Lemma 3.1 at $v$, because $\Gamma$ has it.

Lemma 4.3 provides the inductive step in the proof of our first group-theoretical result.

THEOREM 4.4. *Let $S$ be a finitely generated subgroup of a double $D$ of a LERF group $G$ along a finitely generated subgroup $G_0$, such that the intersection of $S$ with any conjugate of $G_0$ is finitely generated. If $G_0$ is strongly separable in $G$, then $S$ is closed in $\text{PT}(D)$. Hence if $G_0$ is strongly separable in $G$ and if $G_0$ has fgip in $D$, then $D$ is LERF.*

*Proof.* Consider an element $d \in D$ such that $d \notin S$. Remark 3.2 shows that there exists a finite connected subgraph $\Gamma''$ of $\text{Cayley}(D, S)$ representing $S$ and $d$, such that all the vertices of $\Gamma''$ are bichromatic. According to Lemma 3.3 and Lemma 3.6, $\Gamma''$ is contained in a finite subgraph $\Gamma_2$ of $\text{Cayley}(D, S)$ which has properties 4 and 5 of Lemma 3.3, and has property 7 of Lemma 3.5 at any vertex $v$, where $\text{Lab}(\text{Cayley}(D, S), v) \cap G_0$ is of finite index in $G_0$. Then each monochromatic component of $\Gamma_2$ has properties 4', 5' and 7' of Corollary 3.7, and $\Gamma_2$ has properties 2 and 3 of Lemma 3.1; hence $\Gamma_2$ satisfies all the assumptions of Lemma 4.3. Let $U = \{u_3, \cdots, u_n\}$ be the set of all bichromatic vertices of $\Gamma_2$, where $\Gamma_2$ does not have property 7 of Lemma 3.5. Applying Lemma 4.3, we construct a sequence of finite graphs $\Gamma_2, \Gamma_3, \cdots, \Gamma_n$, such that each $\Gamma_i$ is embedded in $\Gamma_{i+1}$, each $\Gamma_i$ has all the properties of Lemma 4.3, and each monochromatic component of $\Gamma_i$ is $G_0$-complete at $u_j \in U$ for $3 \leq j \leq i$. It follows that $\Gamma_n$ has properties 2 and 3 of Lemma 3.1, property 6 of Lemma 3.5, and each monochromatic component of $\Gamma_n$ is $G_0$-complete at any bichromatic vertex of $\Gamma_n$; hence $\Gamma_n$ has property 7 of Lemma 3.5. Therefore Lemma 3.5 and Lemma 3.1 imply that $\Gamma''$ can be embedded in a precover of $D$ with finitely many vertices, and Theorem 4.4 follows from Corollary 2.14.

## 5. Negatively curved groups

A geodesic in a Cayley graph is a shortest path joining two vertices. A group $G$ is $\delta$-negatively curved if any side of any geodesic triangle in Cayley$(G)$ belongs to the $\delta$-neighborhood of the union of the two other sides (see [Gr] and [C-D-P]). We consider only finitely generated negatively curved groups.

Let $\lambda \leq 1$ and $\varepsilon > 0$. A path $p$ in the Cayley graph is a $(\lambda, \varepsilon)$-quasigeodesic if for any subpath $p'$ of $p$ and for any geodesic $\gamma$ with the same endpoints as $p'$, $|\gamma| > \lambda|p'| - \varepsilon$. One of the most important properties of quasigeodesics in negatively curved groups is that for any $\delta$-negatively curved group $G$ and for any pair of numbers $(\lambda, \varepsilon)$, as above, there exists a positive constant $\rho$ which depends only on $(\lambda, \varepsilon)$ and on $\delta$ such that any $(\lambda, \varepsilon)$-quasigeodesic $p$ in Cayley$(G)$ and any geodesic $\gamma$ with the same endpoints as $p$ belong to the $\rho$-neighborhoods of each other (cf. [C-D-P, p. 24].

We use the following property of malnormal quasiconvex subgroups of negatively curved groups proven in [Gi 3]. The original proof of this result for the special case when $G$ is a free group is due to E. Rips ([G-R 2]).

LEMMA 5.1 (the squooshed 4-gon lemma). *Let $G_0$ be a malnormal quasiconvex subgroup of a finitely generated group $G$. Let $\gamma_1 t \gamma_2$ be a path in* Cayley$(G)$ *such that $\gamma_1$ and $\gamma_2$ are geodesics in* Cayley$(G)$, Lab$(\gamma_1) \in G_0$, Lab$(\gamma_2) \in G_0$ *and* Lab$(t) \notin G_0$. *Then for any $L \geq 0$, there exists a positive constant $M(L)$ which depends only on $L$, on $G$ and on $G_0$ such that if $\gamma_1 \subset N_L(\gamma_2)$, then $|\gamma_1| < M(L)$.*

We need the following definitions.

*Definition* 5.2. Let $A = G \underset{G_0=H_0}{*} H$ be as in Definition 2.9. A word $a \equiv a_1 a_2 \cdots a_n \in A$ is in normal form if:

1) $a_i$ lies in one factor of $A$,
2) $a_i$ and $a_{i+1}$ are in different factors of $A$,
3) if $n > 1$, then $a_i \notin G_0$.

Any $a \in A$ has a representative in normal form. If $a \equiv a_1 a_2 \cdots a_n$ is in normal form and $n > 1$, then the Normal Form Theorem ([L-S], p.187) implies that $a$ is not equal to $1_A$.

*Definition* 5.3. Let $p$ be a path in a graph labeled with $X^* \cup Y^*$, and let $p_1 p_2 \cdots p_n$ be its decomposition into maximal monochromatic subpaths. We say that $p$ is in normal form if Lab$(p) \equiv$ Lab$(p_1) \cdots$ Lab$(p_n)$ is in normal form.

Now we will prove our second group-theoretical result.

THEOREM 5.4. *A finitely generated malnormal subgroup $G_0$ of a locally quasiconvex* LERF *negatively curved group $G$ is strongly separable.*

*Proof.* Let $G_0$ be a finitely generated malnormal subgroup of a $\delta$-negatively curved, locally quasiconvex LERF group $G$, let $S_1$ and $S_2$ be finitely generated subgroups of $G$ such that the subgroups $S_1 \cap G_0$ and $S_2 \cap G_0$ are equal and have infinite index in $G_0$, and let $(\Gamma_1, v_1)$ and $(\Gamma_2, v_2)$ be finite subgraphs of covers of $G$ such that $\mathrm{Lab}(\Gamma_1, v_1) = S_1$ and $\mathrm{Lab}(\Gamma_2, v_2) = S_2$ As $S_1$ and $S_2$ are finitely generated subgroups of a locally quasiconvex group $G$, they are quasiconvex in $G$. As the intersection of two quasiconvex subgroups is quasiconvex, the subgroup $S_1 \cap G_0 = S_2 \cap G_0$, which we denote $S_0$, is quasiconvex in $G$, hence it is finitely generated. Choose $K > 0$ such that all the subgroups $G_0, S_1$ and $S_2$ are $K$-quasiconvex. As was mentioned in Lemma 4.1, we consider $(\Gamma_i, v_i)$ as a subgraph of $(\mathrm{Cayley}(G, S_i), S_i \cdot 1), i = 1, 2$. Enlarging $\Gamma_i$, if needed, we can assume that it contains the $K$-neighborhood of $S_i \cdot 1$ in $\mathrm{Cayley}(G, S_i)$. As a locally quasiconvex group has fgip, $\Gamma_i$ is contained in a finite subgraph $\Gamma_i'$ of $\mathrm{Cayley}(G, S_i)$ which has properties $4'$, $5'$ and $7'$ of Corollary 3.7. We will construct embeddings of $(\Gamma_i', v_i)$ in $(\Delta_i, v_i)$ with properties a–d of Definition 4.2. We will use Lemma 6.1 (the ping-pong lemma) proven in Section 6.

Note that as $S_1$ and $S_2$ are $K$-quasiconvex, the constant $C$ described in Lemma 6.1 works for both $S_1$ and $S_2$.

### The construction of $\Delta_i$

As was mentioned at the beginning of this section, there exists a positive constant $\rho$ such that any $(\lambda, \varepsilon)$-quasigeodesic $q$, as in Lemma 6.1, and any geodesic with the same endpoints, as $q$, in $\mathrm{Cayley}(G)$ belong to the $\rho$-neighborhoods of each other.

Let $X_1$ be a finite generating set for $G_0$. Choose a finite generating set $X$ for $G$ such that $X_1$ is a subset of $X$. Using these generating sets, $\mathrm{Cayley}(G_0)$ is a subgraph of $\mathrm{Cayley}(G)$. For any two vertices $g_0'$ and $g_0''$ in $G_0 \subset \mathrm{Cayley}(G)$ consider a geodesic $\gamma$ in $G$ (i.e. a shortest path in $\mathrm{Cayley}(G)$) joining them. Let $\gamma_0$ be a geodesic in $G_0$ (i.e. a shortest path in $\mathrm{Cayley}(G)$ labeled with $X_1^*$) joining $g_0'$ to $g_0''$. As $G$ is locally quasiconvex, $G_0$ is quasiconvex in $G$; hence the embedding of $\mathrm{Cayley}(G_0)$ in $\mathrm{Cayley}(G)$ is a quasi-isometry, and so there exist constants $\lambda_0 \leq 1$ and $\varepsilon_0 > 0$ such that any geodesic $\gamma_0$ in $G_0$, as above, is a $(\lambda_0, \varepsilon_0)$-quasigeodesic in $G$. Thus for any $\gamma$, as above, $|\gamma| > \lambda_0 |\gamma_0| - \varepsilon_0$.

Let $\beta_i$ be the $X_1^*$-component of $S_i \cdot 1$ in $\Gamma_i'$. As the $X_1^*$-component of $S_i \cdot 1$ in $\mathrm{Cayley}(G, S_i)$ is isomorphic to $(\mathrm{Cayley}(G_0, S_0), S_0 \cdot 1)$, we can consider $\beta_1$ and $\beta_2$ as subgraphs of $\mathrm{Cayley}(G_0, S_0)$. Let $B_n$ be the $n$-neighborhood of $S_0 \cdot 1$ in $\mathrm{Cayley}(G_0, S_0)$. There exists a constant $N_0 > 0$ such that for any $n > N_0$, $\beta_1$ and $\beta_2$ are contained in $B_n$, and $\mathrm{Lab}(B_n, S_0 \cdot 1) = S_0$.

Choose a constant $d > 0$ such that $\text{diam}(\Gamma'_i) < d, i = 1, 2$. Let $\lambda_0$, $\varepsilon_0$ and $\rho$ be the constants, defined above, let $L = 2\delta + \rho + d$, let $M(L)$ be as in Lemma 5.1, and let $C$ be as in Lemma 6.1. Choose a constant $N > N_0$ such that $\lambda_0 \cdot N > C + \varepsilon_0$, $\lambda_0 \cdot N > 2(\rho + d) + \varepsilon_0$, and such that $\lambda_0 \cdot N > M(L) + \varepsilon_0$.

As $G_0$ is LERF, and as $S_0$ is finitely generated, there exists an embedding of $(B_N, S_0 \cdot 1)$ in a finite cover $(\alpha_N, v)$ of $G_0$. Let $(\Delta_i, v_i), i = 1, 2$ be the graph constructed from the disjoint union of $(\alpha_N, v)$ and $\Gamma'_i$, with the two copies of $\beta_i$ identified. (Here we use "grafting" again.) Then, by construction, $\Delta_i$ is a finite graph, and the $X_1^*$-component of $v_i$ in $\Delta_i$ is $\alpha_N$, which is a cover of $G_0$. As $\alpha_N$ and $\Gamma'_i$ are isomorphically embedded in $\Delta_i$, we identify $\alpha_N$ and $\Gamma'_i$ with their images in $\Delta_i$.

Denote $G_N = \text{Lab}(\alpha_N, v)$. As $G_0 \cap S_i = S_0 < G_N < G_0$, it follows that $G_N \cap S_1 = G_N \cap S_2 = S_0$. Also, by construction,

$$\text{Lab}(\Delta_i, v_i) = \langle \text{Lab}(\Gamma'_i, v_i), \text{Lab}(\alpha_N, v) \rangle = \langle S_i, G_N \rangle, i = 1, 2.$$

CLAIM 1. *Let $C$ be a constant described in Lemma 6.1. Then all the elements in $G_N$ which are shorter than either $C$ or $2(\rho + d)$, or $M(L)$ (in $G$) belong to $S_0$.*

*Proof.* Consider an element $g \in G_N$ such that $g$ is shorter than $C$ or $2(\rho + d)$, or $M(L)$ (in $G$). Let $\gamma_0$ be a geodesic in $\alpha_N$ which begins at $v$ labeled with $g$. As $\text{Lab}(\alpha_N, v) = G_N$, it follows that $\gamma_0$ ends at $v$. Let $\tilde{\gamma}_0$ be a lift of $\gamma_0$ in $\text{Cayley}(G_0)$. As $\tilde{\gamma}_0$ is a $(\lambda_0, \varepsilon_0)$-quasigeodesic in $\text{Cayley}(G)$, and as the projection map preserves the length of any path, it follows that either $C > |g| > \lambda_0 |\tilde{\gamma}_0| - \varepsilon_0 = \lambda_0 |\gamma_0| - \varepsilon_0$ or, similarly, $2(\rho + d) > |g| > \lambda_0 |\gamma_0| - \varepsilon_0$ or, similarly, $M(L) > |g| > \lambda_0 |\gamma_0| - \varepsilon_0$. In the first case $|\gamma_0| < \frac{C + \varepsilon_0}{\lambda_0} < N$, in the second case $|\gamma_0| < \frac{2(\rho+d)+\varepsilon_0}{\lambda_0} < N$ and in the third case $|\gamma_0| < \frac{M(L)+\varepsilon_0}{\lambda_0} < N$. Hence in all cases $\gamma_0 \subset B_N$, so that $\text{Lab}(\gamma_0) \in S_0 = \text{Lab}(B_N, v)$, proving the claim.

In particular, $G_N$ satisfies the assumptions of Lemma 6.1, and so Lemma 6.1 implies that $\text{Lab}(\Delta_i, v_i) = \langle G_N, S_i \rangle = G_N *_{S_0} S_i$.

We will prove that $\Delta_1$ and $\Delta_2$ have all the properties of Definition 4.2. As these properties can be verified for $\Delta_1$ and $\Delta_2$ separately, we will prove that $\Delta_1$ has them and, to avoid awkward notation, we will write $(\Delta, v)$ for $(\Delta_1, v_1)$, and we will drop the subscript $i = 1, 2$ everywhere else.

Let $u$ and $w$ be vertices in $\Delta$, and let $t$ be a path in $\Delta$ joining $u$ to $w$. There exists a path in $\Delta$ with the same endpoints and the same label as $t$ of the form: $r_1 p r_2$, where $r_1$ and $r_2$ are paths in $\Delta$ joining $u$ to $v$, and $v$ to $w$ respectively, to be specified later, and the path $p$ is a loop beginning at $v$. Hence $\text{Lab}(p) \in \text{Lab}(\Delta, v)$.

Lemma 6.1 implies that there exists a $(\lambda, \varepsilon)$-quasigeodesic $q = q_1 \cdots q_n$ in normal form with $\mathrm{Lab}(p) = \mathrm{Lab}(q)$, such that all $q_i$ are geodesics in $G$.

As $\mathrm{Lab}(r_1)\mathrm{Lab}(q)\mathrm{Lab}(r_2)\mathrm{Lab}(\bar{t}) = 1$, it follows that there exists a closed path $R_1 Q R_2 \bar{T}$ in $\mathrm{Cayley}(G)$ which begins at 1, such that $\mathrm{Lab}(r_i) \equiv \mathrm{Lab}(R_i)$ for $i = 1, 2$, $\mathrm{Lab}(t) \equiv \mathrm{Lab}(T)$ and $\mathrm{Lab}(q) \equiv \mathrm{Lab}(Q)$, where "$\equiv$" means the equality of words. let $Q = Q_1 \cdots Q_n$ be the decomposition of $Q$ with $\mathrm{Lab}(Q_i) \equiv \mathrm{Lab}(q_i)$. Note that if $\mathrm{Lab}(Q_i) \equiv \mathrm{Lab}(q_i) \in S$, then there exists a loop $l$ in $\mathrm{Cayley}(G, S)$ which begins (and ends) at $S \cdot 1$ with $\mathrm{Lab}(l) \equiv \mathrm{Lab}(q_i)$. As $q_i$ is a geodesic in $G$, and as $S$ is $K$-quasiconvex in $G$, the path $l$ belongs to the $K$-neighborhood of $S \cdot 1$ in $\mathrm{Cayley}(G, S)$; thus $l$ belongs to $\Gamma'$. If $\mathrm{Lab}(Q_i) \equiv \mathrm{Lab}(q_i) \in G_N$, then as $\alpha_N$ is a cover of $G_0$, $\alpha_N$ contains a loop $l_\alpha$ which begins at $v$ with $\mathrm{Lab}(l_\alpha) = \mathrm{Lab}(q_i)$. (As $\alpha_N$ is labeled with $X_1^*$, and $q_i$ might not be labeled with $X_1^*$, $\alpha_N$ does not have to contain a loop $l$ with $\mathrm{Lab}(l) \equiv \mathrm{Lab}(q_i)$.) Let $\gamma$ be a geodesic in $G$ with the same endpoints as $Q$. Then $Q$ and $\gamma$ belong to the $\rho$-neighborhood of each other, where $\rho$ is defined after the statement of Lemma 6.1. Hence all $Q_i$ belong to the $\rho$-neighborhood of $\Gamma$.

CLAIM 2. *$\Delta$ has property* a *of Definition 4.2.*

*Proof.* It was mentioned already that $\Delta$ is a finite graph, and that the $X_1^*$-component of $v$ in $\Delta$ is $\alpha_N$, which is a cover of $G_0$. Also, by construction, if an edge $e$ of $\Delta$ does not belong to $\Gamma'$, then $e$ belongs to the $X_1^*$-component of $v$ in $\Delta$. It remains to prove that $\Delta$ can be embedded in a relative Cayley graph of $G$. As $\Delta$ is well-labeled, Lemma 2.16 implies that it is sufficient to show that if vertices $u$ and $w$ in $\Delta$ are joined by a path $t$ with $\mathrm{Lab}(t) = 1$, then $u = w$. Note that if $\mathrm{Lab}(t) = 1$, then $\mathrm{Lab}(p) = \mathrm{Lab}(q) = \mathrm{Lab}(r_1^{-1} r_2^{-1})$, where $p, q$ and $r_i$ are as defined above. Also the path $T$ is a closed loop beginning at 1 in $\mathrm{Cayley}(G)$, so that $R_1 Q R_2$ and $R_1 \gamma R_2$ are closed paths, beginning at 1 in $\mathrm{Cayley}(G)$.

Consider three cases.

1) Both $u$ and $w$ do not belong to the $X_1^*$-component of $v$, so they belong to $\Gamma'$. In this case, let $r_1$ and $r_2$ be the shortest paths in $\Gamma'$ joining $u$ to $v$ and $v$ to $w$ respectively. As $\Gamma'$ has property $5'$, the labels of $r_1$ and $r_2$ are not in $G_0$. As $\mathrm{diam}(\Gamma') < d$, it follows that $r_1$ and $r_2$ are shorter than $d$, and so are $R_1$ and $R_2$. Hence $|\gamma| \leq |R_1| + |R_2| < 2d$. Let $Q = Q_1 \cdots Q_n$ be the decomposition of $Q$, as above. As was mentioned, all $Q_i$ belong to the $\rho$-neighborhood of $\gamma$. As each $Q_i$ is a geodesic in $\mathrm{Cayley}(G)$, it follows that $|Q_i| < 2\rho + |\gamma| < 2\rho + 2d$. If there exists $Q_i$ with $\mathrm{Lab}(Q_i) \in G_N$ then Claim 1 implies that $\mathrm{Lab}(Q_i) \in S_0$. As $\mathrm{Lab}(Q) \equiv \mathrm{Lab}(Q_1) \cdots \mathrm{Lab}(Q_n)$ is in normal form, it follows that $n = 1$ and $\mathrm{Lab}(Q) \equiv \mathrm{Lab}(Q_1) \in S$. Then, as was mentioned above, $\Gamma'$ contains a loop $l$ which begins at $v$ with $\mathrm{Lab}(l) \equiv \mathrm{Lab}(Q_1)$; hence the path $r_1 l r_2$ lies in $\Gamma'$. As $r_1 l r_2$ is labeled with 1 and joins $u$ to $w$, and as $\Gamma'$ is a subgraph of a relative Cayley graph, it follows that $u = w$.

2) If $u$ does not belong to the $X_1^*$-component of $v$, but $w$ does, let $r_1$ be as in case 1, and let $r_2$ be a geodesic in $G_0$ joining $v$ to $w$; hence $\mathrm{Lab}(r_2) \in G_0$. As $\Gamma'$ has property $5'$, any path in $\Gamma'$ joining $u$ to $v$ is labeled by an element not in $G_0$; in particular, $\mathrm{Lab}(r_1) \notin G_0$. Let $R_2'$ be a geodesic in $\mathrm{Cayley}(G)$ with the same endpoints as $R_2$. As $|r_1| = |R_1| < d$, it follows that $\gamma$ belongs to the $(\delta + d)$-neighborhood of $R_2'$; hence $Q$ belongs to the $(\delta + \rho + d)$-neighborhood of $R_2'$. Let $Q = Q_1 \cdots Q_n$ be the decomposition of $Q$, as above. If $Q$ has a subpath $Q_i$ with $i < n$ such that $\mathrm{Lab}(Q_i) \in G_N$ as $(\delta + \rho + d) < L$, Lemma 5.1 implies that $|Q_i| < M(L)$. But then, as in the proof of case 1, $\mathrm{Lab}(Q_i) \in S_0$, hence, as $\mathrm{Lab}(Q) \equiv \mathrm{Lab}(Q_1) \cdots \mathrm{Lab}(Q_n)$ is in normal form, it follows that either $Q = Q_1$ or $Q = Q_1 Q_2$ and $\mathrm{Lab}(Q_2) \in G_N$. If $Q = Q_1 Q_2$ and $\mathrm{Lab}(Q_2) \equiv \mathrm{Lab}(q_2) \in G_N$, then $\mathrm{Lab}(q_1) \equiv \mathrm{Lab}(Q_1) \in S$, $\mathrm{Lab}(r_1)\mathrm{Lab}(q_1)\mathrm{Lab}(q_2)\mathrm{Lab}(r_2) = 1$ and $\mathrm{Lab}(r_1)\mathrm{Lab}(q_1) \in G_0$. As was mentioned above, $\Gamma'$ contains a loop $l$ which begins at $v$ with $\mathrm{Lab}(l) \equiv \mathrm{Lab}(q_1)$, so that $r_1 l$ is a path in $\Gamma'$, which joins $u$ to $v$, and is labeled by an element in $G_0$, contradicting the choice of $u$. If $Q = Q_1$ and $\mathrm{Lab}(Q_1) \in S$, then $\mathrm{Lab}(r_1)\mathrm{Lab}(q_1)\mathrm{Lab}(r_2) = 1$ and $\mathrm{Lab}(r_2) \in G_0$; hence $\mathrm{Lab}(r_1)\mathrm{Lab}(q_1) \in G_0$. Then, as above, $\Gamma'$ contains a path which joins $u$ to $v$ and is labeled by an element in $G_0$, contradicting the choice of $u$. If $Q = Q_1$ and $\mathrm{Lab}(Q) \in G_N$, then $\mathrm{Lab}(r_1) \in G_0$, contradicting the choice of $r_1$. Hence case 2 cannot occur.

3) Both $u$ and $w$ belong to the $X_1^*$-component of $v$; hence $u$ and $w$ are in $\alpha_N$. Let $r_1$ and $r_2$ be geodesics in $\alpha_N$ joining $u$ to $v$ and $v$ to $w$ respectively. Then $\mathrm{Lab}(r_1) \equiv \mathrm{Lab}(R_1) \in G_0$, and $\mathrm{Lab}(r_2) \equiv \mathrm{Lab}(R_2) \in G_0$. Now $\mathrm{Lab}(\gamma) \in G_0$. Let $Q = Q_1 \cdots Q_n$ be the decomposition of $Q$, as above. If $Q$ has a subpath $Q_i$ with $1 < i < n$ such that $\mathrm{Lab}(Q_i) \in G_N$ then, as $Q_i$ belongs to the $\rho$-neighborhood of $\gamma$, and $\rho < L$, Lemma 5.1 implies that $|Q_i| < M(L)$. But then Claim 1 implies that $\mathrm{Lab}(Q_i) \in S_0$; hence, as $\mathrm{Lab}(Q) \equiv \mathrm{Lab}(Q_1) \cdots \mathrm{Lab}(Q_n)$ is in normal form, it follows that either $Q = Q_1$ or $Q = Q_1 Q_2$. If $Q = Q_1 Q_2$, assume, without loss of generality, that $\mathrm{Lab}(Q_2) \in G_N$, and $\mathrm{Lab}(Q_1) \in S$. Then $\mathrm{Lab}(R_1)\mathrm{Lab}(Q_1)\mathrm{Lab}(Q_2)\mathrm{Lab}(R_2) = 1$. As $\mathrm{Lab}(Q_2)\mathrm{Lab}(R_2)\mathrm{Lab}(R_1) \in G_0$, it follows that $\mathrm{Lab}(Q_1) \in G_0$. Then $\mathrm{Lab}(Q_1) \in G_0 \cap S = S_0$, contradicting the definition of the normal form.

If $Q = Q_1$ and $\mathrm{Lab}(Q_1) \in S$, then as $\mathrm{Lab}(R_2)\mathrm{Lab}(R_1) \in G_0$, it follows that $\mathrm{Lab}(Q_1) \in G_0 \cap S = S_0 < G_N$. Hence, $Q$ should have the form $Q = Q_1$ with $\mathrm{Lab}(Q_1) \in G_N$. Then, as was mentioned above, $\alpha_N$ contains a loop $l_\alpha$ which begins at $v$ with $\operatorname{Lab}(l_\alpha) = \operatorname{Lab}(Q_1)$, and so $r_1 l_\alpha r_2$ is a path in $\alpha_N$. As $r_1 l_\alpha r_2$ is labeled with 1 and joins $u$ to $w$, and as $\alpha_N$ is a cover, it follows that $u = w$.

CLAIM 3. $\operatorname{Lab}(\Delta, v) \cap G_0 = G_N$, hence $\operatorname{Lab}(\Delta_1, v_1) \cap G_0 = \operatorname{Lab}(\Delta_2, v_2) \cap G_0$, therefore $\Delta_1$ and $\Delta_2$ have property b of Definition 4.2.

*Proof.* Claim 2 demonstrated that $\Delta$ can be considered as a subgraph of $\operatorname{Cayley}(G, G_N *_{S_0} S)$. As $\operatorname{Lab}(\Delta, v) = G_N *_{S_0} S$, it follows that $\operatorname{Lab}(\Delta, v) \cap G_0$ contains $G_N$. Consider a loop $a$ in $\Delta$ which begins at $v$ with $\operatorname{Lab}(a) \in G_0$. As $\alpha_N$ is a subgraph of $\Delta$, and $\alpha_N$ is $G_0$-complete at $v$, there exists a path $a'$ in $\alpha_N$ which begins at $v$ with $\operatorname{Lab}(a) = \operatorname{Lab}(a')$. As $\Delta$ is a subgraph of a relative Cayley graph, it follows that $a$ and $a'$ should have the same terminal vertex, namely $v$. Thus $a'$ is a loop in $\alpha_N$ beginning at $v$ and $\operatorname{Lab}(a) = \operatorname{Lab}(a') \in \operatorname{Lab}(\alpha_N, v) = G_N$. Therefore $\operatorname{Lab}(\Delta, v) \cap G_0$ is contained in $G_N$, proving Claim 3.

CLAIM 4. $\Delta$ has property d of Definition 4.2.

*Proof.* Let $u$ and $w$ be vertices in the image of $\Gamma'$ in $\operatorname{Cayley}(G, G_N *_{S_0} S)$, and let $t$ be a path in $\operatorname{Cayley}(G, G_N *_{S_0} S)$ joining $u$ and $w$, such that $\operatorname{Lab}(t) \in G_0$. Let $T'$ be a geodesic in $\operatorname{Cayley}(G)$ with the same endpoints as $T$, where $T$ is as above. Consider three cases listed in the proof of Claim 2.

1) In this case $|r_i| = |R_i| < d$; hence $Q$ belongs to the $(2\delta + \rho + d)$-neighborhood of $T'$. Let $Q = Q_1 \cdots Q_n$ be the decomposition of $Q$, as above. If there exists $Q_i$ with $\operatorname{Lab}(Q_i) \in G_N$ then, as $(2\delta + \rho + d) < L$, Lemma 5.1 implies that $|Q_i| < M(L)$. Then, as in the proof of Claim 2, it follows that $Q = Q_1$, and $\operatorname{Lab}(Q_1) \in S$. Then, as was mentioned above, $\Gamma'$ contains a loop $l$ which begins at $v$ with $\operatorname{Lab}(l) \equiv \operatorname{Lab}(Q_1)$; hence the path $r_1 l r_2$ lies in $\Gamma'$. As $r_1 l r_2$ joins $u$ and $w$ and as $\operatorname{Lab}(r_1 l r_2) = \operatorname{Lab}(t) \in G_0$, property d holds.

2) In this case we can assume that $w = v$ and $r_2$ is a trivial path. As $|R_1| < d$, it follows that $Q$ belongs to the $(\rho + \delta + d)$-neighborhood of $T'$. Let $Q = Q_1 \cdots Q_n$ be the decomposition of $Q$, as above. If there exists $i < n$ such that $\operatorname{Lab}(Q_i) \in G_0$, then as in the proof of Claim 1, $\operatorname{Lab}(Q_i) \in S_0$. Hence either $Q = Q_1$ or $Q = Q_1 Q_2$ and $\operatorname{Lab}(Q_2) \in G_N$. But then, as in the proof of Claim 2, $\Gamma'$ should contain a path labeled with $G_0$ which joins $u$ to $v$, contradicting the choice of $u$. Hence case 2 cannot occur.

3) Let $r_1$ and $r_2$ be as in case 3 of Claim 2. Then the path $r_1 r_2$ joins $u$ to $w$, and its label is in $G_0$, and therefore property d holds.

CLAIM 5. $\Delta$ has property c of Definition 4.2.

*Proof.* Let $u$ be a vertex in the image of $\Gamma'$ in $\Delta$ which does not belong to the $X_1^*$-component of $S \cdot 1$ in Cayley$(G, S)$. Let $r_1$ be a shortest path in $\Gamma'$ joining $u$ to $v$. Let $t$ be a loop in $\Delta$ beginning at $u$ with $\mathrm{Lab}(t) \in G_0$. Let $r_2 = \bar{r}_1$, let $R_1, R_2, Q$ and $T$ be as in the proof of case 1 of Claim 4 with $u = w$. Then, as in the proof of Claim 4, it follows that $Q = Q_1$, and $\mathrm{Lab}(Q_1) \in S$; hence $\Gamma'$ contains a loop $l$ which begins at $v$ with $\mathrm{Lab}(l) \equiv \mathrm{Lab}(Q_1)$; so the path $r_1 l r_2 = r_1 l \bar{r}_1$ lies in $\Gamma'$. Then $r_1 l \bar{r}_1$ is a loop in $\Gamma'$ which begins at $u$ with $\mathrm{Lab}(r_1 l \bar{r}_1) = \mathrm{Lab}(t) \in G_0$. As $\Gamma'$ has property $4'$, it follows that $\Delta$ has property c of Definition 4.2.

## 6. The Ping-Pong Lemma

Let $\lambda \leq 1, \mu > 0$ and $\varepsilon > 0$. A path $p$ is a local $(\lambda, \varepsilon, \mu)$-quasigeodesic if for any subpath $p'$ of $p$ which is shorter than $\mu$ and for any geodesic $\gamma$ with the same endpoints as $p', |\gamma| > \lambda |p'| - \varepsilon$ (cf. [C-D-P, p. 24].

Theorem 1.4 (p. 25) of [C-D-P] (see also [Gr, p. 187]) states that for any $\lambda' \leq 1$ and for any $\varepsilon' > 0$ there exist constants $(\mu, \lambda, \varepsilon)$ which depend only on $(\lambda', \varepsilon')$ and $\delta$, such that any local $(\lambda', \varepsilon', \mu)$-quasigeodesic in $G$ is a global $(\lambda, \varepsilon)$-quasigeodesic in $G$.

LEMMA 6.1. *Let $S$ and $G_0$ be $K$-quasiconvex subgroups of a $\delta$-negatively curved group $G$, and let $S_0 = S \cap G_0$. If $G_0$ is malnormal in $G$, then there exists a constant $C > 0$ which depends only on $G, \delta$ and $K$, such that for any subgroup $G_N$ of $G_0$ with $G_0 \cap S = G_N \cap S = S_0$, if all the elements in $G_N$ which are shorter than $C$ (in $G$) belong to $S_0$, then the following hold:*

1) $\langle G_N, S \rangle = G_N *_{S_0} S$, *where $\langle G_N, S \rangle$ denotes the minimal subgroup of $G$, containing $G_N$ and $S$.*

2) *There exist constants $\lambda \leq 1$ and $\varepsilon > 0$ such that for any element in $\langle G_N, S \rangle$ there exists a $(\lambda, \varepsilon)$-quasigeodesic representative (in $G$) $q = q_1 q_2 \cdots q_m$ in normal form, where all $q_j$ are geodesics in $G$.*

*Proof.* Let $G_N$ be a subgroup of $G_0$ such that $G_N \cap S = G_0 \cap S = S_0$. Let $l$ be an element of $\langle G_N, S \rangle$ such that $l \notin S_0$. Then $l$ can be written as a product $l = g_1 s_1 \cdots s_{m-1} g_m$, where $g_i \in G_N, s_i \in S, g_i$ and $s_i$ do not belong to $S_0$, $g_i$ and $s_i$ are geodesics in $G$, $g_1$ is a shortest representative of the coset $g_1 S_0$, $g_m$ is a shortest representative of the coset $S_0 g_m$, and for $1 < i < m, g_i$ is a shortest representative of the double coset $S_0 g_i S_0$. (The elements $g_1$ or $g_m$ might be trivial.) Let $p$ be the path in Cayley$(G)$ beginning at 1 which has the form $p = p_1 q_1 \cdots q_{m-1} p_m$, where $\mathrm{Lab}(p_i) \equiv g_i$ and $\mathrm{Lab}(q_i) \equiv s_i$.

Let $A$ be the number of words in $G$ which are shorter than $2K + \delta$ and let $M(2\delta)$ be as in Lemma 5.1. Let $\lambda' = 1/6$ and $\varepsilon' = 4K \cdot A + \delta + M(2\delta)$.

As mentioned above, there exist constants $(\mu, \lambda, \varepsilon)$ which depend only on $\lambda'$, $\varepsilon'$ and on $\delta$ such that any local $(\lambda', \varepsilon', \mu)$-quasigeodesic in $G$ is a global $(\lambda, \varepsilon)$-quasigeodesic in $G$. (In this case $(\mu, \lambda, \varepsilon)$ depend only on the constants $1/6, K, A, M(2\delta)$ and $\delta$.)

Let $C = \max(\mu, \frac{\varepsilon}{\lambda})$. We claim that if all elements in $G_N$ which are shorter than $C$ belong to $S_0$, then any path $p$, as above, is a $(\lambda, \varepsilon)$-quasigeodesic in $G$. Indeed, it is enough to show that $p$ is a local $(1/6, (4K \cdot A + \delta + M(2\delta)), \mu)$-quasigeodesic in $G$.

As $g_i \notin S_0$, it follows that $|p_i| > C$. As $\mu < C$, any subpath $t$ of $p$ with $|t| < \mu$ has a (unique) decomposition $t_1 t_2 t_3$, where $t_1$ and $t_3$ are subpaths of some $p_i$ and $p_{i+1}$, and $t_2$ is a subpath of $q_i$ (some of $t_i$ might be empty). Let $t_4$ be a geodesic in $G$ connecting the endpoints of $t$. By definition, $p$ is a local $(1/6, (4K \cdot A + \delta + M(2\delta)), \mu)$-quasigeodesic in $G$ if and only if $|t_4| \geq \frac{|t|}{6} - (M(2\delta) + \delta + 4K \cdot A)$.

If 2 out of the 3 subpaths $t_i$ are empty, then the remaining one, say $t_2$, is a geodesic; therefore $|t| = |t_2| = |t_4|$, and $t_4$ satisfies the above inequality. If at least 2 out of the 3 subpaths $t_i$ are nonempty, then $t_2$ is nonempty. Considering $l^{-1}$ instead of $l$, if needed, we can assume that $t_1$ is nonempty.

If $|t_2| > \frac{2|t|}{3}$, then $|t_1| + |t_3| \leq \frac{|t|}{3}$, so that $|t_4| \geq |t_2| - (|t_1| + |t_3|) \geq \frac{2|t|}{3} - \frac{|t|}{3} = \frac{|t|}{3} > \frac{|t|}{6} - (M(2\delta) + \delta + 4K \cdot A)$.

If $|t_2| \leq \frac{2|t|}{3}$, assume, without loss of generality, that $|t_1| \geq |t_3|$; then $|t_1| > \frac{|t|}{6}$. As $t_1 t_2 t_3 t_4$ is a geodesic 4-gon in a $\delta$-negatively curved group $G$, there exists a decomposition $t_1 = t'_2 t'_3 t'_4$ such that $t'_2$ belongs to the $\delta$-neighborhood of $t_2$, $t'_3$ belongs to the $2\delta$-neighborhood of $t_3$ and $t'_4$ belongs to the $\delta$-neighborhood of $t_4$. According to Lemma 5.1, $|t'_3| < M(2\delta)$ and according to Lemma 6.2 (below), $|t'_2| \leq 4K \cdot A$. But then $|t_4| + \delta \geq |t'_4| = |t_1| - |t'_2| - |t'_3| \geq |t_1| - 4K \cdot A - M(2\delta) \geq \frac{|t|}{6} - 4K \cdot A - M(2\delta)$. Hence $|t_4| \geq \frac{|t|}{6} - (M(2\delta) + \delta + 4K \cdot A)$, so the path $p$ is a $(\lambda, \varepsilon)$-quasigeodesic in $G$.

As $C \geq \frac{\varepsilon}{\lambda}$, it follows that if a $(\lambda, \varepsilon)$-quasigeodesic $p$ in $\text{Cayley}(G)$ with $\iota(p) = 1$ is longer than $C$, then $|1, \text{Lab}(p)| \geq \lambda|p| - \varepsilon > \lambda C - \varepsilon > 0$. Hence $\text{Lab}(p)$ is not equal to 1. As any element $l \in \langle G_N, S \rangle$ which is not in $S_0$ has a representative $\text{Lab}(p)$ in $G$, as above, with $|p| > C$, it follows that $l$ is not equal to 1; hence $\langle G_N, S \rangle = G_N *_{S_0} S$.

LEMMA 6.2. *When the notation of the proof of Lemma* 6.1 *is used,* $|t'_2| \leq A \cdot 4K$.

*Proof.* To simplify notation, we drop the subscript $i$ on the paths $p_i$ and $q_i$ and on their labels, so $t_1$ is a subpath of $p$, $t_2$ is a subpath of $q$, $\text{Lab}(p) = g$ and $\text{Lab}(q) = s$. As $G_N < G_0$, we consider $g$ as an element of $G_0$. Without loss of generality, assume that $q$ begins at 1 (so it ends at $s$); then $p$ begins at $g^{-1}$ and ends at 1. As $G_0$ and $S$ are $K$-quasiconvex in $G$, any vertex $v_i$ on $p$ is in the $K$-neighborhood of $G_0$, and any vertex $w_i$ on $q$ is in the $K$-neighborhood of $S$. Hence we can find vertices $v_1$ and $v_2$ in $t_2'$, $w_1$ and $w_2$ in $t_2$, $g'$ and $g''$ in $G_0$, and $s'$ and $s''$ in $S$ such that $|v_i, w_i| < \delta, |v_1, (g')^{(-1)}| < K, |v_2, (g'')^{(-1)}| < K, |w_1, s'| < K$ and $|w_2, s''| < K$. We choose $v_1$ between 1 and $v_2$, and we choose $w_1$ between 1 and $w_2$. Let $\gamma'$ be a geodesic in Cayley$(G)$ joining $(g')^{-1}$ to $s'$, and let $\gamma''$ be a geodesic in Cayley$(G)$ joining $(g'')^{-1}$ to $s''$. Then Lab$(\gamma') = g's'$, Lab$(\gamma'') = g''s''$, $|\gamma'| < 2K + \delta$ and $|\gamma''| < 2K + \delta$.

Assume that $|t_2'| > A \cdot 4K$. Then we can find vertices, as above, which, in addition, satisfy: $|v_2, v_1| > 4K$ and Lab$(\gamma') =$ Lab$(\gamma'')$; hence $g's' = g''s''$. It follows that $(g'')^{(-1)}g' = s''(s')^{(-1)}$, so both products are in $S_0$. As $g$ is a shortest element in the double coset $S_0 g S_0$, it follows that $|g| \leq |g(g'')^{(-1)}g'|$. Let $r$ be a geodesic joining $(g'')^{(-1)}$ to $v_2$, let $b'$ be a subpath of $p$ joining $v_2$ to 1 and let $b''$ be a subpath of $p$ joining $g^{-1}$ to $v_2$. Then $|g| = |p| = |b'| + |b''|$, and $|g(g'')^{(-1)}g'| \leq |b''| + |r| + |g'|$; hence $|b'| + |b''| \leq |b''| + |r| + |g'|$, so that $|b'| + |r| \leq 2|r| + |g'|$. As $|g''| \leq |b'| + |r|$, and as $|r| \leq K$, it follows that $|g''| \leq 2K + |g'|$.

As $|v_2, v_1| > 4K$, the triangle inequality implies that $|g''| = |(g'')^{-1}| \geq |b'| - |r| = |1, v_1| + |v_1, v_2| - |r| \geq |1, v_1| + 4K - K = |1, v_1| + 3K$.

Let $a$ be a geodesic joining $(g')^{-1}$ to $v_1$. As $|a| < K$, the triangle inequality implies that $|g'| = |(g')^{-1}| \leq |1, v_1| + |a| < |1, v_1| + K$. Hence, $|g''| > |g'| + 2K$, a contradiction. Therefore $|t_2'| \leq A \cdot 4K$.

*Acknowledgement.* The author would like to thank her friends for their support.

CALIFORNIA INSTITUTE OF TECHNOLOGY, PASADENA, CA
*Current address*: A & H CONSULTANTS, ANN ARBOR, MI
*E-mail address*: ritagtk@math.lsa.umich.edu

## REFERENCES

[A-D]   R. B. J. T. ALLENBY and D. DONIZ, A free product of finitely generated nilpotent groups amalgamating a cycle that is not subgroup separable, Proc. A.M.S. **124** (1996), 1003–1005.

[A-G]   R. B. J. T. ALLENBY and R. J. GREGORAC, On locally extended residually finite groups, *Conference on Group Theory*, Lecture Notes in Math. **319**, Springer-Verlag, New York, 1973, 9–17.

[B-F]   M. BESTVINA and M. FEIGHN, A combination theorem for negatively curved groups, J. Differential Geom. **35** (1992), 85–101.

[Bo]    F. BONAHON, Bouts des variétés hyperboliques de dimension 3, Ann. of Math. **124** (1986), 71–158.

[B-B-S] A. M. BRUNNER, R. G. BURNS, and D. SOLITAR, The subgroup separability of free products of two free groups with cyclic amalgamation, Contemp. Math. **33**, A.M.S., Providence, R. I., 1984, 90–115.

[B-K-S] R. G. Burns, A. Karrass, and D. Solitar, A note on groups with separable finitely generated subgroups, Bull. Austral. Math. Soc. **36** (1987), 153–160.
[Ca] R. Canary, A covering theorem for hyperbolic 3-manifolds and its applications, Topology **35** (1996), 751–778.
[Car] L. Carrol, *Through the Looking-Glass and what Alice Found There*, Progress Publishers, Moscow, 1966.
[C-D-P] M. Coornaert, T. Delzant, and A. Papadopoulos, *Géométrie et Théorie des Groupes*, Lecture Notes in Math. **1441**, Springer-Verlag, New York, 1990.
[Ge] S. Gersten, Coherence in doubled groups, Comm. Algebra **9** (1981), 1893–1900.
[Gi 1] R. Gitik, Graphs, LERF groups and the topology of 3-manifolds, Ph. D. Thesis, Hebrew University (1990).
[Gi 2] ———, Graphs and separability properties of groups, J. of Algebra **188** (1997), 125–143.
[Gi 3] ———, On quasiconvex subgroups of negatively curved groups, J. Pure Appl. Algebra **119** (1997), 155–169.
[Gi 4] ———, Ping-pong on negatively groups, J. of Algebra **217** (1999), 65–72.
[Gi 5] ———, On the profinite topology on doubles of groups, J. of Algebra **219** (1999), 80–86.
[Gi 6] ———, On the combination theorem for negatively curved groups, Internat. J. Algebra Comput. **6** (1996), 751–760.
[G-R 1] R. Gitik and E. Rips, A Necessary condition for $A \underset{a=b}{*} B$ to be LERF, Israel J. Math. **73** (1991), 123–125.
[G-R 2] ———, On separability properties of groups, Internat. J. Algebra Comput. **5** (1995), 703–717.
[Gre] L. Greenberg, Discrete groups of motions, Canad. J. Math. **12** (1960), 415–426.
[Gro] M. Gromov, Hyperbolic groups, *Essays in Group Theory*, 75–263, MSRI series **88**, S. M. Gersten, ed., Springer-Verlag, New York, 1987.
[G-T] J. L. Gross and T. W. Tucker, *Topological Graph Theory*, Wiley-Interscience, New York, 1987.
[Hall] M. Hall, Jr., Coset representations in free groups, Trans. A.M.S. **67** (1949), 421–432.
[Hi] G. Higman, A finitely related group with an isomorphic proper factor group, J. London Math. Soc. **26** (1951), 59–61.
[Jaco] W. Jaco, *Lectures on Three-Manifold Topology*, A.M.S. Publications **43**, Providence, R.I., 1980.
[L-N] D. D. Long and G. A. Niblo, Subgroup separability and 3-manifold groups, Math. Z. **207** (1991), 209–215.
[L-S] R. C. Lyndon and P. E. Schupp, *Combinatorial Group Theory*, Springer-Verlag, New York, 1977.
[Ma] B. Maskit, *Kleinian Groups*, Springer-Verlag, New York, 1987.
[Ni] G. Niblo, Separability properties of free groups and surface groups, J. Pure Appl. Algebra **78** (1992), 77–84.
[Rips] E. Rips, On a double of a free group, Israel J. of Math. **96** (1996), 523–525.
[R-W] J. H. Rubinstein and S. Wang, $\pi_1$-injective surfaces in graph manifolds, Comment. Math. Helv. **73** (1998), 499–515.
[Sco 1] P. Scott, Subgroups of surface groups are almost gometric, J. London Math. Soc. **17** (1978), 555–565.
[Sco 2] ———, Correction to Subgroups of surface groups are almost geometric, J. London Math. Soc. **32** (1985), 217–220.
[Sta] J. R. Stallings, Topology of finite graphs, Invent. Math. **71** (1983), 551–565.

[Thu]  W. P. THURSTON, *Three-dimensional manifolds, Kleinian groups and hyperbolic geometry*, Bull. A.M.S. **6** (1982), 357–381.
[We]   B. A. F. WEHRFRITZ, *The Residual finiteness of some generalised free products*, J. London Math. Soc. **24** (1981), 123–126.